\documentstyle{article}
\renewcommand{\baselinestretch}{1.3}

\newtheorem {th}{Theorem}[section]
\newtheorem {lem}[th]{Lemma}
\newtheorem {pr}[th]{Proposition}
\newtheorem {cor}[th]{Corollary}

\def\tensor{\otimes}
\def\VS{{{\cal V}_S}}
\def\VG{{{\cal V}_G}}
\def\VZ{{{\cal V}_{\Z^d}}}
\def\V0{{{\cal V}_{0}}}
\def\VT{T^d}
\def\xy{{\vec {xy}}}
\def\zw{{\vec {zw}}}
\def\uv{{\vec {uv}}}

\def\Cox{\hfill \Box}

\def\deq{\, {\stackrel {def} {=}}}

\def\dconv{\, {\stackrel {{\cal D}} {\rightarrow}}}

\def\dd{\delta}

\def\ee{\epsilon}
\def\E{{\bf{E}}}
\def\P{{\bf{P}}}
\def\pois{{\cal P}_1}
\def\N{\hbox{I\kern-.2em\hbox{N}}}
\def\R{\hbox{I\kern-.2em\hbox{R}}}
\def\Z{\hbox{Z\kern-.4em\hbox{Z}}}
\def\CC{\hbox{C\kern-.45em\hbox{\rule{.01in}{1.4ex}}\hspace{.45em}}}

\def\C{{\cal{C}}}

\def\|{\, | \, }

\def\Tree{{\bf T}}

\begin{document}

\begin{titlepage}
\begin{center}
{\large \bf LOCAL CHARACTERISTICS, ENTROPY AND LIMIT THEOREMS FOR
SPANNING TREES AND DOMINO TILINGS VIA TRANSFER-IMPEDANCES}
\end{center}
Running Head: LOCAL BEHAVIOR OF SPANNING TREES 
\vspace{5ex}
\begin{flushright}
Robert Burton \footnote{This research supported in part by Grant number
AFOSR 90-2015} \\
Robin Pemantle \footnote{This research supported by a National 
Science Foundation postdoctoral fellowship .  } \footnote{Now
at the University of Wisconsin-Madison}\\
Oregon State University \\

\end{flushright}

\vfill

{\bf ABSTRACT:} \break
Let $G$ be a finite graph or an infinite graph on which $\Z^d$ acts
with finite fundamental domain.  If $G$ is finite, let $\Tree$ be a
random spanning tree chosen uniformly from all spanning trees of $G$; 
if $G$ is infinite, methods from \cite{Pe} show that this still
makes sense, producing a random {\em essential spanning
forest} of $G$.  A method for calculating 
local characteristics (i.e. finite-dimensional marginals) of $\Tree$
from the {\em transfer-impedance} matrix is presented.  This differs
from the classical matrix-tree theorem in that only small pieces of the
matrix ($n$-dimensional minors) are needed to compute small
($n$-dimensional) marginals.  Calculation of the matrix entries
relies on the calculation of the Green's function for $G$, which is not
a local calculation.  However, it is shown how the calculation
of the Green's function may be reduced to a finite computation in the case
when $G$ is an infinite graph admitting a $Z^d$-action with finite
quotient.  The same computation also gives the entropy of the law 
of $\Tree$.  

These results are applied to the problem of tiling certain lattices by
dominos -- the so-called {\em dimer problem}.  Another application of
these results is to prove modified versions of conjectures of Aldous 
\cite{Al2} on the limiting distribution of degrees of a vertex
and on the local structure near a vertex
of a uniform random spanning tree in a lattice whose dimension
is going to infinity.  Included is a generalization of moments
to tree-valued random variables and criteria for these generalized
moments to determine a distribution.

\vfill



\end{titlepage}

\section{Introduction}
 
We discuss spanning trees and domino tilings (perfect matchings)
of periodic lattices.  To define these terms, let $S = 
\{ 1 , \ldots , k \}$ be a generic $k$-element set and
let $G$ be a graph whose vertex set is $\Z^d \times k$, 
i.e. its vertices are all pairs $(x,i)$ where $x= (x_1 , 
\ldots , x_d)$ is a vector of integers of length $d$ and
$i$ is an integer between $1$ and $k$.  We will usually allow $G$ to
denote the vertex set $\Z^d \times S$, since this causes
no ambiguity.  We say that $G$
is {\em periodic} if its edge set is invariant under the natural
$\Z^d$-action; in other words we require that $(x,i)$ is
connected to $(y,j)$ (written $(x,i) \sim (y,j)$) if and only
if $(0,i) \sim (y-x,j)$.  Assume throughout that $G$ is connected
and locally finite.  By periodicity there is a maximal degree
$D$ of vertices of $G$.  Graphs that we consider may have
parallel edges, in other words more than one edge may connect
the same two vertices.  It is convenient to add self-edges -- edges
connecting a vertex to itself -- until all vertices have degree
$D$; such a graph is called {\em $D$-regular}.  Adding self-edges
does not alter any of the problems we address, so 
we assume throughout that all vertices of $G$ have degree $D$.
It will also be convenient to assume that simple random walk
on $G$ is aperiodic.  Since this is true whenever $G$ has a self-edge,
we assume the presence of at least one self-edge.

A spanning tree of any graph is subcollection of the edges
having no loops, but such that every pair of vertices is connected
within the subcollection.  A loopless subgraph that is not 
necessarily connected is called a forest, and a forest in
which every vertex is connected to infinitely many others is
called an essential spanning forest.  It is shown in \cite{Pe} that
the uniform measures on spanning trees of a cube of finite 
size $n$ in the integer lattice $\Z^d$ converge weakly as 
$n \rightarrow \infty$ to a measure $\mu_{\Z^d}$ on essential spanning
forests of $\Z^d$.  This measure chooses a spanning tree
with probability one if $d \leq 4$ and with probability
zero if $d \geq 5$.  The first purpose of the present work
is to show how the finite dimensional marginals and the entropy
of the limiting measure $\mu_{\Z^d}$ may be effectively computed.  Since 
the computations may be carried out in the more general setting
of periodic lattices, and since some of these (e.g. the
hexagonal lattice in the plane) seem as interesting as
$\Z^d$ from the point of view of physical modelling, we
treat the problem in this generality.  The methods of \cite{Pe}
extend to the case of arbitrary periodic graphs to show that
$\mu_G$ chooses a spanning tree with probability one if $d \leq 4$ 
and zero if $d \geq 5$.  We will not re-prove this result in the
more general setting, since that would involve a completely 
staight-forward but lengthy redevelopment of the theory of 
loop-erased random walks \cite{La2} for periodic graphs.

A domino tiling of a graph is a partition of the vertices into
sets of size two, each set containing two adjacent vertices.
Domino tilings on $\Z^2$ have been studied \cite{Kas} and the
exponential growth rates of the number of tilings of large regions with
various boundary conditions have been calculated.  The growth rates
for domino tilings are different for different boundary conditions
\cite{Kas,TF,Elk}.  The second purpose
of this work is to exhibit the domino tiling of maximal entropy 
for each periodic lattice in a special class (that includes $\Z^2$), 
and to compute the entropy.  We exploit a general version 
of a known connection between domino tilings and spanning forests, 
so that this follows more or less immediately from the 
results on spanning forests.  The correspondence also gives a way
to calculate probabilities of various contours arising in 
a uniform random domino tiling.  The first few of these are 
calculated by Fisher \cite{Fi1,Fi2} using Pfaffians.

The method we use to calculate the f.d.m.'s of $\mu_G$ is to calculate
the Green's function for $G$ and then to write the f.d.m.'s 
as determinants of the {\em transfer-impedance matrix}, which is
a matrix of differences of the Green's function.  The main result on 
transfer impedance matrices is stated and proved in Section 4
(Theorem~\ref{transfer}) in the general setting of periodic 
lattices.  Since the result is interesting in itself and is for
the rest of the paper {\em sine qua non}, we state here a simplified 
version for finite graphs.  
\begin{th} \label{simplified}
Let $G$ be a finite $D$-regular graph.  Fix an arbitrary orientation
of the edges of $G$ and for edges $e = \xy$ and $f$ of $G$ define
$H(e,f)$ to be the expected signed number of transits of $f$ by
a random walk started at $x$ and stopped when it hits $y$ (this
can be written as a difference of Green's functions).  For edges
$e_1 , \ldots , e_k$ let $M(e_1, \ldots , e_k)$ denote the matrix
whose $i,j$-entry is $H(e_i , e_j)$.  If $\Tree$ is a random spanning
tree, uniformly distributed among all spanning trees of $G$, then
$$\P (e_1 , \ldots , e_k \in \Tree ) = \det M(e_1 , \ldots , e_k ) . $$
$\Cox$
\end{th}
The idea of a transfer-impedance matrix is not new,
the terminology being taken from \cite{We}.  We have not, however,
been able to find the key result (Theorem~\ref{transfer}) on
determinants of the transfer-impedance matrix stated anywhere.
Furthermore, extending results about transfer-impedances
from finite graphs to infinite graphs is not immediate (at least
when simple random walk on the infinite graph  is transient)
and requires an argument based on triviality of the Poisson boundary.
For these reasons, we include a derivation of all results on
transfer-impedances from scratch.  

The rest of the paper is organized as follows.  The next section
contains the notation used in the rest of the paper and a
derivation of the 
Green's function for a periodic lattice.  Section 3 contains 
lemmas, such as a discrete Harnack's principle, about simple random
walks on periodic lattices.  Rather than providing detailed
proofs, we include in an appendix the outline of a standard proof
 for the case $G = \Z^d$ and 
indicate the necessary modifications for arbitrary periodic
lattices.  The connection between simple random walks and spanning
trees is documented in \cite{Pe}; the main result that
will be used from there is that $\P (e \in \Tree)$ is determined
by certain hitting probabilities, but the reader desiring more details
may also consult \cite{Al2} or \cite{Br}.  
Section 4 uses these lemmas to show that the 
Green's function is the unique limit of Green's functions
on finite subgraphs and that it in fact determines the
f.d.m.'s for $\mu_G$ via determinants of the transfer-impedance
matrix.  Section 5 considers two examples.  The first is
the case $G = \Z^2$, which is special because the Green's
function for general periodic lattices is given by a definite
integral which is only explicitly evaluable when $G= \Z^2$.
The second is the high dimensional limit of $G = \Z^d , d
\rightarrow \infty$, which converges in a sense to be defined later
to a critical Galton-Watson Poisson(1) branching process,
in accordance with a conjecture of Aldous \cite{Al1}.  Section 6
calculates the entropy of $\mu_G$.  Section 7 discusses the
connection between spanning trees of a lattice and domino
tilings of the join of a lattice and its dual.  From this
follows a determination of the topological entropy for domino tilings
of graphs that are joins of a periodic lattice and its dual.  
It is also possible from this to exhibit f.d.m.'s of the 
maximal entropy domino tiling in a few special cases.  

Certain characterizations of the Green's function and Harnack principles
are required that are essentially adaptations of known results on $\Z^d$
to arbitrary periodic lattices.  These adaptations are treated briefly
in the appendix.  Also given in the appendix are criteria for
determining limits of probability distributions on trees from knowledge
of certain functionals which act as generalized moments.  

\section{Notation and a Green's Function}

Let $G$ be a periodic lattice with the assumptions of
connectedness, $D$-regularity and there being at least one loop,
as in the previous section.  Let $\mu_G$ denote the
weak limit as $n \rightarrow \infty$ of the uniform
measures on spanning trees of the induced subgraph on
$G$ with vertices $\{ (x,i) :\, \parallel x \parallel_\infty \,
\leq n \}$.  The arguments in \cite{Pe} show that this limit exists 
because the probability of the elementary event of a finite
set of edges all being in the tree is always decreasing in
$n$ and because the measures of these elementary 
events determine $\mu_G$.  It follows from Corollary~\ref{converges},
the random walk construction of spanning trees and
the first equality of Lemma~\ref{gives answer} that,
as in the case where $G = \Z^d$, the limit can be taken
independent of the boundary conditions, i.e. the limit 
for induced subgraphs is the same as for tori.  For an edge
$e$, we often write $\P (e \in \Tree )$ for $\mu \{ \Tree :
e \in \Tree \}$.  

For any finite set of edges $e_1 , \ldots
, e_k$ that form no loop among them there is a graph $G / e_1 , 
\ldots , e_k$ called the contraction of $G$ by $e_1 , \ldots , e_k$.
Its vertices are the vertices of $G$ modulo the equivalence relation
of being connected by edges in $\{ e_1 , \ldots , e_k \}$.  Let the
projection from vertices in $G$ to vertices of $G / e_1 , \ldots ,
e_k$ be called $\pi$.  Then the edges of $G / e_1 , \ldots , e_k$
are precisely one edge connecting $\pi (x)$ to $\pi (y)$ for each edge 
connecting $x$ to $y$ in $G$.  If $W$ is a subset of the vertices
of $G$, the {\em induced subgraph} of $G$ on $W$ is the graph with
vertex set $W$ and an edge $\xy$ for each edge $\xy \in G$ with $x,y
\in W$.  For use in Section 7 we include the dual notion to contraction,
namely deletion.  If $e_1 , \ldots , e_k$ are edges of connected graph
$G$ whose removal does not disconnect $G$ then the deletion of $G$ by 
$e_1 , \ldots , e_k$, denoted $G-e_1 , \ldots , e_k$
is simply $G$ with $e_1 , \ldots e_k$ removed.  The salient point is
that deletion and contraction 
by different edges commute, so that $G/e_1 , \ldots , e_k - e_1',
\ldots , e_l'$ is well-defined when $e_i \neq e_j'$ for all $i,j$.

Let $SRW_x^G$ denote simple random walk on $G$ starting from $x$.
More precisely, a path in $G$ is a function $f$ from the nonnegative 
integers to the vertices of $G$ such that $f(i+1)$ is always adjacent
to $f(i)$; $SRW_x^G$ makes all possible intial segments of
a given length equally likely.  Write $\P (SRW_x^G (i) = y)$
for the probability that a simple random walk on $G$ started from
$x$ is at $y$ at time $i$.  Either $G$ or $x$ may be suppressed 
in the notation when no ambiguity arises.  For a subset $B$
of the vertices of $G$, let $\partial B$ denote the boundary of
$B$, namely those vertices $x \in B$ that have neighbors in $B^c$.
For any $x \notin B$, let $\tau_x^B = \inf \{ j : SRW_x (j) \in B \}$
denote the (possibly infinite) hitting time of $B$ from $x$.  

Define the vector spaces $\VS , \VZ$ and $\VG$ to be the set of
complex-valued functions on $S$, $\Z^d$ and $\Z^d \times S$ respectively,
with pointwise addition, scalar multiplication, and the topology of
pointwise convergence.  For $u \in \VS$ and $g \in \VZ$ let
$u \tensor g$ denote the $f \in \VG$ for which $f(x,i) = g(x) u(i)$.
Let $\VT$ denote the $d$-dimensional torus $\R^d / \Z^d$, written as
$d$-tuples of elements of $(-1/2,1/2]$.  For $\alpha \in \VT$ and $x \in \Z^d$,
the inner product $\alpha \cdot x \deq \sum_{i=1}^d \alpha_i x_i$ is
a well-defined element of $\R / \Z$; for fixed $\alpha \in \VT$,
let $\xi^\alpha \in \VZ$  be defined by $\xi^\alpha (x) = 
e^{2 \pi i \alpha \cdot x}$.  Define the adjacency operator $A : \VG
\rightarrow \VG$ by $(Af) (x,i) = D^{-1} \sum_{(y,j) \sim (x,i)} f(y,j)$.
Here, and in all subsequent such summations, the element $(y,j)$
is to be counted as many times as there are edges from $(x,i)$ to $(y,j)$.
A function $f \in \VG$ is called harmonic if $Af = f$ and harmonic at
$(x,i)$ if $(Af) (x,i) = f(x,i)$.  We define a family of $k$ by $k$ adjacency 
matrices $\{ R^x : x \in \Z^d \}$ by letting $R^x (i,j)$ equal one if $(0,i) 
\sim (x,j)$ and zero otherwise.  Observe that $R^x = 0$ for all but 
finitely many $x \in \Z^d$ and that $R^{-x} = (R^x)^T$.
For $\alpha \in \VT$, define the $k$ by $k$
matrix $Q (\alpha)$ by the essentially finite sum $D^{-1} \sum_{x \in \Z^d} 
e^{2 \pi i \alpha \cdot x} R^x$.  

Now we begin building Green's functions for simple random walk on $G$.
For a transient walk, the Green's function $H(x,y)$ can be defined as
the expected number of visits to $y$ starting from $x$.  This is
symmetric, and harmonic in each argument except on the diagonal.  Here,
we construct a function $g^f$ that is harmonic except at a finite number
of points $x$, at which $(I-A)g^f (x)$ is equal to some specified $f$.  
Later, a uniqueness
theorem will show that when $f = \dd_x$ for some $x \in G$, then
$g^f$ specializes to $H(x,\cdot )$.  In dimension
two, simple random walk is recurrent.  In this case, although the
Green's function may still be defined classically by subtracting
the expected number of visits from $x$ to itself, the integral defining $g^f$
will blow up for $f = \dd_x$.  It will, however, be finite for $f = \dd_x -
\dd_y$.

When $G$ is just $\Z^d$ with the usual nearest-neighbor edges, the
eigenfunctions in $\VG$ for the adjacency operator are just the
functions $\xi^\alpha$.  The first lemma uses these to construct 
eigenfunctions for $A$.

\begin{lem} \label{eigenfunctions}
Suppose $u \in \VS$ satisfies $Q (\alpha) (u) = \lambda u$ for some
real $\lambda$.  Then $A(u \tensor \xi^\alpha) = \lambda u \tensor \xi^\alpha$.
\end{lem}

\noindent{Proof:} 
\begin{eqnarray*}
A (u \tensor \xi^\alpha) (x,i) & = & {1 \over D} \sum_{(y,j) \sim (x,i)}
    (u \tensor \xi^\alpha) (y,j) \\[2ex]
& = & {1 \over D} \sum_{z \in \Z^d ,1 \leq j \leq k} R^z (i,j) u_j 
    e^{2 \pi i \alpha \cdot (x+z)} \\[2ex]
& = & e^{2 \pi i \alpha \cdot x} \sum_{1 \leq j \leq k} ({1 \over D}
    \sum_{z \in \Z^d} e^{2 \pi i \alpha \cdot z} R^z) (i,j) u_j \\[2ex]
& = & e^{2 \pi i \alpha \cdot x} \sum_{1 \leq j \leq k} Q(\alpha)
    (i,j) u_j \\[2ex]
& = & e^{2 \pi i \alpha \cdot x} \lambda u_i \\[2ex]
& = & \lambda (u \tensor \xi^\alpha) (x,i)   \Cox
\end{eqnarray*}

This lemma tells us how to invert $(I-A)$ on elements of $\VG$ of the
form $u \tensor \xi^\alpha$ where $u$ is an eigenvector of $Q(\alpha)$
with eigenvalue not equal to one.
By representing general elements of $\VG$ as integrals of eigenfunctions
we can then invert $(I-A)$ on these integrals since $(I-A)^{-1}$ 
commutes with the integral, at least when absolute integrability
conditions are satisfied.  A preliminary observation is that for
any $u \in \VS$,
$$ \int u \tensor \xi^\alpha  \; d\alpha = u \tensor \dd_0 ,$$
where $d\alpha$ is the usual Haar measure on $\VT$.  To see this,
note that the integrand is bounded in magnitude by $|u|$, so the
integral makes sense and integrating pointwise gives
$$\int u \tensor \xi^\alpha (x,i) \; d\alpha = u_i \int e^{2 \pi i \alpha
    \cdot x} d\alpha = u_i \dd_0(x) .$$
Inverting $(I-A)$ on elements of $\VG$ with finitely many nonzero coordinates
is easily reduced to inverting $(I-A)$ on things of the form $u \tensor
\dd_0$ and their translates by elements of $\Z^d$.  
The above representation shows that these are integrals of
$u \tensor \xi^\alpha$ which are sums of eigenfunctions (the eigenfunctions
given by letting u be an eigenvector of $Q(\alpha))$.  So the above 
representation solves the problem as long as the inverted integrand, 
$(1 - \lambda)^{-1} u \tensor \xi^\alpha$ remains integrable.  With this 
in mind, observe that $Q(\alpha)$ is Hermitian for each $\alpha$.  It
is therefore diagonalizable with real eigenvalues and
has a unitary basis of eigenvectors.  Let $\{v (\alpha , i) : 
\alpha \in \VT , 1 \leq i \leq k \}$ denote a measurable selection
of ordered eigenbasis for $Q(\alpha )$ and let $\lambda (\alpha , i)$
be the eigenvalue corresponding to $v(\alpha , i)$.  Now for $u \in \VS$,
let $c^u (\alpha , i) = \ll u , v(\alpha , i) \gg$ 
denote the coefficients of $u$ in the chosen eigenbasis, in other words
$$ \sum_i c^u (\alpha , i) v(\alpha , i) = u $$
for each $\alpha \in \VT$.  Then we have the following theorem.

\begin{th} \label{exhibit Green}
For $u \in \VS$, let 
\begin{equation} \label{gu}
g^u = \int \sum_{i=1}^k Re \left \{ c^u (\alpha , i) (1 - \lambda 
    (\alpha , i))^{-1} v(\alpha , i) \tensor \xi^\alpha
    \right \} \; d\alpha 
\end{equation}
and for $j \leq d$, define
\begin{equation} \label{guj}
g^{u,j} = \int \sum_{i=1}^k Re \left \{ c^u (\alpha , i) (1 - \lambda 
    (\alpha , i))^{-1} (1 - e^{-2 \pi i \alpha_j}) v(\alpha , i) \tensor 
    \xi^\alpha \right \} \; d\alpha .
\end{equation}
The integrals are meant pointwise, i.e. as defining $g^u (x)$ and
$g^{u,j} (x)$ for each $x \in G$.  Then
\begin{quote}
$(i)$ The integrand in~(\ref{guj}) is always integrable and the integrand
in~(\ref{gu}) is integrable when $d \geq 3$ or when $d=1$ or $2$ and
$\sum_i u_i = 0$. \\

$(ii)$ $(I-A) g^u = u \tensor \dd_0$ and $(I-A) g^{u,j} = u \tensor
(\dd_0 - \dd_{e_j})$ whenever the integrals exist, $e_j$ being the
$j^{th}$ standard basis vector in $\Z^d$. \\

$(iii)$ $g^u$ and $g^{u,j}$ are bounded whenever the defining 
integrals exist.
\end{quote}
\end{th}

\noindent{Remark:} When $d \geq 2$, it is not necessary to take the
real part in~(\ref{gu}) and~(\ref{guj}) since the imaginary part
integrates to zero.  When $d=1$ however, the imaginary part fails to
be integrable.

For $f,g \in \VG$, say $f$ is a translate of $g$ if $f(x,i) = g(x+a,i)$
for some $a \in \Z^d$.  If $f \in \VG$ has finitely many nonzero 
coordinates, then it can be represented as a finite sum of translates 
of elements of the form $u \tensor \dd_0$.  Further, if the sum of
$f(x,i)$ is zero, then $f$ can be represented as the sum of translates
of elements $u \tensor \dd_0$ for which $\sum u_j = 0$ together with
translates of elements $u \tensor (\dd_0 - \dd_{e_j})$.  Thus the following 
is an immediate corollary.

\begin{cor} \label{general Green}
Let $f \in \VG$ have finitely many nonzero coordinates.  If $d \geq 3$
or $d=1$ or $2$ and $\sum_{x \in G} f(x) = 0$, then the previous two theorems
can be used to construct a bounded solution $g$ to $(I-A) g = f$.  $\Cox$
\end{cor}

The proof of Theorem~\ref{exhibit Green}
depends on the following lemma which bounds the eigenvalues
of $Q(\alpha)$ away from $1$ in terms of $|\alpha|$ in order to
get the necessary integrability results.

\begin{lem} \label{inverse square}
There is a constant $K = K(G)$ for which $\max_i |\lambda (\alpha ,i)|
\leq 1 - K |\alpha|^2$, where for specificity we take $|\alpha|
= \max_j |\alpha_j|$.
\end{lem}

\noindent{Proof:}  The eigenvalues of a matrix are continuous functions
of its entries \cite{Ka}, so it suffices to show that this is true in a 
neighborhood of $0$ and to show that $\lambda (\alpha , i) \neq 1$
for $\alpha \neq 0$.  For the first of these, it suffices to find
for each $s \leq d$ a constant $K_s$ for which the eigenvalues of $Q(\alpha)$ 
are bounded in magnitude by $1 - K_s |\alpha_s|^2$.  So fix an $s \leq d$.

Begin with a description of the entries of $Q(\alpha)^r$.  The  
quantity $Q(\alpha)(i,j)$ is the sum over edges of $G$ connecting $(0,i)$
to $(x,j)$ for $x \in \Z^d$ of complex numbers of modulus $1/D$.
Furthermore, as $i$ or $j$ varies with the other fixed, there are
precisely $D$ of these paths.  It follows that $Q(\alpha)^r (i,j)$
is the sum over paths of length $r$ connecting $(0,i)$ to some
$(x,j)$ of complex numbers of modulus $D^{-r}$ and that there
are $D^r$ of these contributions in every row and every column.

Suppose we can find an $r = r(s)$ such that for every $i \leq k$ there
is a path of length $r$ from $(0,i)$ to $(e_s,i)$ where $e_s$ is
the $s^{th}$ standard basis vector.   Since there is a self-edge
at every vertex, there is perforce a path of length $r$ from
$(0,i)$ to $(0,i)$.  These two paths represent summands in the
above decompositon of $Q(\alpha)^r (i,i)$ whose arguments differ
by $\alpha_s$.  By the law of cosines, the sum of these two 
terms has magnitude $D^{-r}(2+2\cos (\alpha_s))^{1/2} \leq
2 D^{-r} (1-c \alpha_s^2)$ for any $c < 1/8$ and $\alpha$ in 
an appropriate neighborhood of zero.  Adding in the rest of
the terms in the $i^{th}$ row of $Q(\alpha)^r$ and using the triangle
inequality shows that the sum of the magnitudes of the entries in the
$i^{th}$ row is
most $1 - 2cD^{-r} \alpha_s^2$ for $\alpha_s$ in a neighborhood of
zero.  The usual Perron-Frobenius argument then shows that no 
eigenvalue has greater modulus than the maximal row sum of moduli.  
This implies that in an appropriate neighborhood of zero,
no eigenvalue of $Q(\alpha)$ has modulus greater
than $(1-2cD^{-r}\alpha_s)^{1/r} \leq 1-K_s \alpha_s^2$
where $K_s = 2cD^{-r(s)}/r(s)$, which is the bound we wanted.

Finding such an $r$ is easy.  Since $G$ is connected there is
for each for each $i \leq k$ a path of some length $l_i$ from
$(0,i)$ to $(e_s,i)$.  These can be extended to paths of
any greater length by including some self-edges, so $r$
can be chosen as the maximum over $i$ of $l_i$.  All that remains
is to show that the only time $Q(\alpha)$ has an eigenvalue of 1 
is when $\alpha = 0$.  This is essentially the same argument.  
Picking $r' \geq \max_s r(s)$, the row sums of the moduli of
of the entries of $Q(\alpha)^{r'}$ are strictly less than
one and the Perron-Frobenius argument shows that no eigenvalue
has modulus one or greater.    $\Cox$

\noindent{Proof} of Theorem~\ref{exhibit Green}:  First we establish
when the integrands in~(\ref{gu}) and~(\ref{guj}) are integrable.
For each $\alpha$  the vectors $v(\alpha , i)$ form a unitary basis,
hence the coefficients $c^u (\alpha , i)$ are bounded in magnitude
by $|u|$.  Since $v(\alpha , i)$, $\xi^\alpha$ and $1 - e^{-2 \pi 
\alpha_j}$ all have unit modulus, integrability of 
$|(1 - \lambda (\alpha , i))^{-1}|$ is certainly sufficient to
imply integrability of~(\ref{gu}) and~(\ref{guj}).  For 
$d \geq 3$, this now follows immediately from $|(1 - \lambda 
(\alpha , i))^{-1}| \leq K |\alpha^{-2}|$.

For $d = 2$, $|\alpha|^{-2}$ is not integrable, so we will
find another factor in the integrands that is $O(|\alpha|)$.
By aperiodicity of SRW on
$G$, we also have that the projected SRW on $S$ is aperiodic, and
therefore that the eigenvalue of $1$ when $\alpha = 0$ is simple
with eigenvector $v(0,1) = D^{-1/2} (1, \ldots , 1)$.  (Here we assume
without loss of generality that the eigenvectors have been numbered
so that for $\alpha$ in a neighborhood of zero, $v(\alpha,1)$ 
is an eigenvector whose eigenvalue has maximum modulus.)  The
assumption $\sum_i u_i = 0$ in~(\ref{gu}) then implies that 
$c^u (0,1) = 0$.  By analyticity of the eigenbasis
with respect to the entries of the matrix (at least away
from multiple eigenvalues) \cite{Ka}, $c^u (\alpha , 1) = 
O(|\alpha|)$.  Thus $|c^u (\alpha , i) (1 - \lambda (\alpha , i))^{-1}|
\leq K' |\alpha|^{-1}$ which implies integrabilty 
of~(\ref{gu}) when $d=2$.  Similarly, $(1-e^{2 \pi i 
\alpha_j}) = O(|\alpha|)$, which implies integrability of~(\ref{guj})
when $d=2$.

When $d=1$, we will show that the real parts of
\begin{equation} \label{eq3}
[c^u (\alpha , 1) v(\alpha , 1) \tensor \xi^\alpha ] (x,i)
\end{equation}
for $\sum_i u_i = 0$ and
\begin{equation} \label{eq4}
[c^u (\alpha , 1) v(\alpha , 1) (1 - e^{-2 \pi i \alpha}) \tensor 
    \xi^\alpha ] (x,i)
\end{equation}
for any $u$ are both $O(|\alpha|^2)$ as $\alpha \rightarrow 0$
for fixed $(x,i) \in G$.  Clearly, this is enough to imply integrability
of~(\ref{gu}) and~(\ref{guj}).  Observe that~(\ref{eq3}) and~(\ref{eq4})
are both zero when $\alpha = 0$.  By analyticity of $v(\alpha , 1)$
at zero, it suffices to show that derivatives of~(\ref{eq3}) 
and~(\ref{eq4}) with respect to $\alpha$ at zero are purely
imaginary.  Taking~(\ref{eq4}) first, we have
\begin{eqnarray*}
&& {d \over d\alpha} \left ( c^u (\alpha , 1) v(\alpha , 1) (1 -
    e^{-2 \pi i \alpha }) \tensor \xi^\alpha (x,i) \right ) |_{\alpha
    = 0} \\[2ex]
& = & 2 \pi i c^u (0,1) v(0,1) \tensor \xi^0 (x,i) = 2 \pi i \ll u , v(0,1) \gg 
    v(0,1)_i \in \sqrt{-1} \R ,
\end{eqnarray*}
since the factor of $1 - e^{-2 \pi i \alpha}$ kills all the other
terms in the derivative when $\alpha = 0$.  For~(\ref{eq3}), we get
\begin{eqnarray*}
&& {d \over d\alpha} \left ( c^u (\alpha , 1) v(\alpha , 1) 
    \tensor \xi^\alpha (x,i) \right ) |_{\alpha = 0} \\[2ex]
& = & [v(0,1) \tensor \xi^0] (x,i) {d \over d\alpha} c^u (\alpha , 1)
    |_{\alpha = 0} \\[2ex]
& = & v(0,1)_i \ll u , (d / d\alpha)|_{\alpha = 0} v(\alpha , 1) \gg
\end{eqnarray*}
so it suffices to show that $v(\alpha , 1)$ has imaginary derivative
at zero. 

Observe first that $Q(\alpha)$ has imaginary derivative at $\alpha = 0$
since the entries of $Q(\alpha)$ are all sums of $e^{2 \pi i \alpha x}$
for various $x \in \Z$.  Call this imaginary derivative $R$.
Secondly, observe that the derivative of $\lambda (\alpha , 1)$ vanishes
at $\alpha = 0$ since $\lambda (\alpha , 1)$ is real and attains its
maximum at $\alpha = 0$.  Then letting $w$ denote the derivative
of $v(\alpha , 1)$ at $\alpha = 0$,
$$[Q + \ee R + O(\ee^2)] (v + \ee w + O(\ee^2)) = (1 + O(\ee^2))
    (v + \ee w + O(\ee^2)) $$
from which it follows that
$$ Rv = (I - Q) w + O(\ee) . $$
Letting $\ee \rightarrow 0$ gives $Rv = (I-Q)w$.  Since $R$ is imaginary
and $I,Q$ and $v$ are real, it follows that $w$ is imaginary.  This
shows that the real part of~(\ref{eq3}) is $O(|\alpha|^2)$ and
completes the proof of $(i)$.

The above argument actually also establishes boundedness of $g^u$
and $g^{u,j}$ when $d \geq 2$, but we give a different probabilistic
argument since it is necessary to do so anyway for the case $d=1$.
Let $s \in (0,1)$ be a real parameter and consider the functions 
$g_s^u$ and $g_s^{u,j}$ gotten by replacing $\lambda (\alpha , i)$
by $s \lambda (\alpha , i)$ in~(\ref{gu}) and~(\ref{guj}).  The
integrands are {\em a fortiori} absolutely integrable, being
bounded in magnitude by $|1 - s|^{-1}$ times a possible factor
of $2$ for the $1 - e^{-2 \pi i \alpha_j}$.  Thus in fact
$g_s^u \in l^\infty (G)$.  Now taking the real part is no longer
necessary, since the imaginary part is an odd function of
$\alpha$ and must integrate to zero.  We have then
\begin{eqnarray} \label{killedRW1}
g_s^u & = & \int \sum_{i=1}^k \sum_{n=0}^\infty s^n \lambda 
    (\alpha , i)^n c^u (\alpha , i) v(\alpha , i)
    \tensor \xi^\alpha \; d\alpha \\[2ex]
& = & \sum_{n=0}^\infty s^n \int \sum_{i=1}^k \lambda 
    (\alpha , i)^n c^u (\alpha , i) v(\alpha , i)
    \tensor \xi^\alpha \; d\alpha \nonumber \\[2ex]
& = & \sum_{n=0}^\infty s^n \int \sum_{i=1}^k c^u (\alpha , i)
    A^n (v(\alpha , i) \tensor \xi^\alpha) \; 
    d\alpha \nonumber \\[2ex]
& = & \sum_{n=0}^\infty s^n A^n \left ( \int \sum_{i=1}^k c^u 
    (\alpha , i) (v(\alpha , i) \tensor \xi^\alpha) \right ) \; 
    d\alpha \nonumber \\[2ex]
& = & \sum_{n=0}^\infty s^n A^n (u \tensor \dd_0)  . \nonumber
\end{eqnarray}
The reason $A$ may be commuted with the sum and integral is that
$A f (x)$ is a finite linear combination of terms $f(y)$ for $y \in G$,
and each of these terms is integrable.  In a similar manner, we get
\begin{equation} \label{killedRW2}
g_s^{u,j} = \sum_{n=0}^\infty s^n A^n (u \tensor (\dd_0 - \dd_{e_j})) .
\end{equation}
Since $A^n$ gives the transition probabilities for an $n$-step
simple random walk on $G$, this means that
$$g_s^u (x) = \sum_{n=0}^\infty \E u \tensor \dd_0 (X_n) $$
where $X_n$ is a $SRW_x^G$ killed with probability $1-s$ at each step.
Similarly,
$$g_s^{u,j} (x) = \sum_{n=0}^\infty \E u \tensor (\dd_0 - \dd_{e_j}) (X_n). $$
It follows from this that $g_s^u , g_s^{u,j} \rightarrow 0$ as
$x \rightarrow \infty$ for fixed $s,u,j$.  From the forward equation
for the random walk (or by direct calculation from~(\ref{killedRW1})
and~(\ref{killedRW2})), $g_s^u = sAg_s^u + u
\tensor \dd_0$ and $g_s^{u,j} = sAg_s^{u,j} + u \tensor (\dd_0 - \dd_{e_j})$,
whence it follows that $g_s^u$ (resp. $g_s^{u,j}$) cannot have a maximum
or minimum except on the support of $u \tensor \dd_0$ (resp. $u \tensor
(\dd_0 - \dd_{e_j})$).  Since  $u \tensor \dd_0$ (resp. $u \tensor
(\dd_0 - \dd_{e_j})$) has finite support, say $W \subset G$, this
implies that for all $y \in G$,
\begin{equation} \label{g_s bounded}
\min_{x \in W} g_s^u(x) \leq g_s^u (y) \leq \max_{x \in W} g_s^u(x) ,
\end{equation}
and similarly for $g_s^{u,j}$.
The proof of integrability of~(\ref{gu}) and~(\ref{guj}) shows that 
$g_s^u \rightarrow g^u$ and $g_s^{u,j} \rightarrow g^{u,j}$ as
$s \rightarrow 1$.  Taking the limit of~(\ref{g_s bounded}) gives
$$\min_{x \in W} g^u(x) \leq g^u (y) \leq \max_{x \in W} g^u(x) $$
and similarly for $g^{u,j}$, establishing $(iii)$.

Finally to show $(ii)$, we have from~(\ref{killedRW2}) that
$$(I-sA) g_s^u = u \tensor \dd_0 .$$
Since $g_s^u \rightarrow g^u$ pointwise as $s \rightarrow 1$
and since $A g_s^u (x)$ is a finite sum of values $g_s^u (y)$,
the limit of the LHS as $s \rightarrow 1$ exists and is equal to
$(I-A) g^u$.  A similar argument for $g^{u,j}$ completes the
proof of $(ii)$ and of the theorem.   $\Cox$

\section{Simple Random Walk on G}

Define the cube $B_n$ of size $n$ in $G$ to be those vertices $(x,i)$ 
for which $|x_j| \leq n : 1 \leq j \leq d$.  For $x \in B_n$, define
the hitting distribution on the boundary, $\nu_{x}^{B_n}$, to be 
the law of $SRW_{x} (\tau_{x}^{B_n})$.  Thus for example, if
$x \in \partial B_n$, then $\nu_{x}^{B_n} = \dd_{x}$.
For $x \notin B_n$, define the hitting distribution on the boundary 
to be the same, but conditioned on the SRW hitting the boundary; thus
$\nu_{x}^{B_n} (C) = \P(SRW_{x} (\tau_{x}^{B_n}) \in C) 
/ \P (\tau_{x}^{B_n} < \infty)$.  For $n < m$ and $x \in \partial
(B_m^c)$, let $\rho_x^{B_n B_m}$ be the hitting distribution on 
$\partial B_n$ of a $SRW_x$ conditioned never to return to $B_m^c$;
for $m < n$ and $x \in \partial B_m$, let $\rho_x^{B_n B_m}$ be
the hitting distribution on $\partial B_n$ of $SRW_x$ conditioned 
not to return to $B_m$.  The following lemma is an adaptation of
the discrete Harnack inequalities on $\Z^d$ for general periodic
lattices.  The proof is merely an adaptation of the proof for $\Z^d$
and will be sketched in the appendix.

\begin{lem}[Harnack principles] \label{Harnack}
Let $G$ be a periodic graph satisfying the assumptions of the
first section.  
Fix a positive integer $n$.  Then as $m \rightarrow \infty$,
\begin{eqnarray*}
(i) && \max_{x,y \in B_n , z \in \partial B_m} \nu_x^{B_m} (\{z\}) / 
    \nu_y^{B_m} (\{z\}) \rightarrow 1 \\[2ex]
(ii) && \max_{x,y \in \partial B_n , z \in \partial B_m} \rho_x^{B_m B_n} 
    (\{z\}) / \rho_y^{B_m B_n} (\{z\}) \rightarrow 1 \\[2ex]
(iii) &&\max_{z \in \partial B_n , x,y \in B_m^c} \nu_x^{B_n} (\{z\}) / 
    \nu_y^{B_n} (\{z\}) \rightarrow 1 \\[2ex]
(iv) &&\max_{z \in \partial B_n , x,y \in \partial (B_m^c)} 
    \rho_x^{B_n B_m} (\{z\}) / \rho_y^{B_n B_m} (\{z\}) \rightarrow 1 \\[2ex]
(v) && (i) - (iv) \mbox{ hold when $G$ is replaced by a finite contraction 
    of } G
\end{eqnarray*}
$\Cox$
\end{lem}

\begin{cor} \label{const}
Bounded harmonic functions on finite contractions of
periodic graphs are constant.
\end{cor}

\noindent{Proof:}  For any vertex $x$, Let $X_0 , X_1 , \ldots$
be a simple random walk starting from $x$.  If $g$ is harmonic
then $\{ g(X_i) \}$ is a martingale, and if $g$ is bounded
and $x \in B_m$ then optional stopping gives $g(x) = 
\E g(X(\tau_x^{B_m})) = \int g(z) d \nu_x^{B_m} (z)$.  
By $(i)$ of the previous lemma, $\nu_x^{B_m} (z) = (1+0(1))
\nu_y^{B_m} (z)$ as $m \rightarrow \infty$, hence $g(x)
= g(y)$ and $g$ is constant.   $\Cox$

\begin{cor} \label{unique}
Let $G$ be a finite contraction of a periodic graph.  Then for $f \in \VG$
there is, up to an additive constant at most one bounded solution
$g$ to $(I-A)g = f$.
\end{cor}

\noindent{Proof:}  If $g_1$ and $g_2$ are two solutions then $g_1 
- g_2$ is a bounded harmonic function.   $\Cox$

\begin{cor} \label{converges}
Let $G$ be a finite contraction of a periodic graph and let
$G_n$ be the induced subgraph on $B_n$.  For $x,y,z \in G$ with
$x \sim y$ define $h(x,y,z,n) = \P (SRW_z^{G_n} (\tau_z^y - 1) =x)$
to be the probability that $SRW_z^{G_n}$ first hits $y$ by coming
from $x$ (with $h(x,y,y,n) \deq 0$).  
Then $\lim_{n \rightarrow \infty} h(x,y,z,n)$ exists
for all $x,y,z$.
\end{cor}

\noindent{Proof:}  Fix $L$ such that $x,y \in B_L$.  For $w \in B_L$
define
$$\phi_1 (x,y,w) = \P (\tau_w^y < \infty \mbox{ and } SRW_w
    (\tau_w^y - 1) = x)$$
and 
$$\phi_2 (x,y,w) = \P (\tau_w^y < \infty \mbox{ and } SRW_w
    (\tau_w^y - 1) \neq x) . $$
For $i = 1,2,$ let $\phi_i (x,y,w,n)$ be $\phi_i (x,y,w)$ with
the clause $\tau_w^y < \infty$ replaced by $\tau_w^y < \tau_w^{B_n^c}$
and observe that $\phi_i (x,y,w,n) \rightarrow \phi_i (x,y,w)$
as $n \rightarrow \infty$.
From the Harnack lemma we know that $\rho_z^{B_L B_m}$ approaches
a limiting measure $\rho$ on $\partial B_L$ as $m \rightarrow \infty$.  
We claim that
\begin{eqnarray} \label{infinity}
h(x,y,z,n) & \rightarrow & \P (\tau_z^y < \infty \mbox{ and } SRW_z
    (\tau_z^y - 1) = x) \\
& + & \P (\tau_z^y = \infty) \int \phi_1 (w) /
    (\phi_1 (w) + \phi_2 (w)) \, d\rho (w) . \nonumber
\end{eqnarray}
To see this, write $h(x,y,z,n)$ as $\P (SRW_z^G (\tau_z^y - 1) = x
\mbox{ and } \tau_z^y < \tau_z^{G_n}) + \P (\tau_z^y > \tau_z^{G_n})
\P (SRW_z^{G_n} (\tau_z^y - 1) = x \| \tau_z^y > \tau_z^{G_n})$.
The first of these terms is clearly converging to the first term
in~(\ref{infinity}), while the first factor of the second is
converging to $\P (\tau_z^y = \infty)$; the second factor is a mixture
over $u \in \partial B_n$ of $\P (SRW_u^{G_n} (\tau_u^y - 1) = x)$,
so it suffices to show that this is converging to the integral
in~(\ref{infinity}) uniformly in $u$ and $n \rightarrow \infty$.  
Consider the sequence of times $\tau_1 , \sigma_1 , \tau_2 , \sigma_2,
\ldots$ where $\tau_1$ is the first time that $SRW_u$ hits $B_L$, 
$\sigma_1$ is the nest time it hits $\partial B_n$, $\tau_2$ is the next
time it hits $B_L$, and so forth.  The first hitting time $\tau_u^y$
of $y$ must satisfy $\tau_i \leq \tau_u^y < \sigma_i$ for some $i$.
Now write
$$ \P (SRW_u^{G_n} (\tau_u^y - 1) = x) =$$
$$ { \sum_{i = 1}^\infty \P (SRW_u^{G_n} (\tau_u^y - 1) = x ; \tau_i \leq
     \tau_u^y < \sigma_i) \over
     \sum_{i = 1}^\infty \P (SRW_u^{G_n} (\tau_u^y - 1) \neq x ; \tau_i \leq
     \tau_u^y < \sigma_i) + \P (SRW_u^{G_n} (\tau_u^y - 1) = x ; \tau_i \leq
     \tau_u^y < \sigma_i) } .$$
The sum in the denominator is of course $1$, but the point of writing
it this way is to illustrate that for each $i$ the ratio is approximately
the integral in~(\ref{infinity}).  More precisely, for fixed $i$ 
the Markov property gives that
$\P (SRW_u^{G_n} (\tau_u^y - 1) = x ; \tau_i \leq \tau_u^y < \sigma_i)$
is equal to $\P (\tau_u^y > \sigma_{i-1})$ times a mixture over $v \in
\partial B_n$ (corresponding to the last exit from $\partial B_n$ before
$\tau_i$) of $\int \phi_1 (x,y,w,n) \, d\rho_v^{B_L B_n}$.  Similarly,
$\P (SRW_u^{G_n} (\tau_u^y - 1) \neq x ; \tau_i \leq \tau_u^y < \sigma_i)$
is equal to $\P (\tau_u^y > \sigma_{i-1})$ times a mixture over $v \in
\partial B_n$ of $\int \phi_2 (x,y,w,n) \, d\rho_v^{B_L B_n}$.  Since
$\phi_i (x,y,w,n) \rightarrow \phi_i (x,y,w)$ and $\rho_v^{B_L B_n}
= (1+o(1)) \rho$ as $n \rightarrow \infty$, this shows that 
the ratio of the numerator to the denominator in the sum is $(1+o(1))$
times the integral in~(\ref{infinity}) and proves the corollary.    
$\Cox$

\section{Transfer Impedance}

Before going into the definition of transfer impedance, it is worth
pausing to remark that the functions $g^u$ constructed in
Theorem~\ref{exhibit Green} really are versions of the Green's 
function.  This is not essential to any of the arguments below, 
so the proofs are relegated to the appendix.  Define the usual
Green's function $H(\cdot,\cdot)$ on pairs of vertices of a 
periodic graph $G$ by 
$$H(x,y) = \sum_{n=0}^\infty \P (SRW_x (n) = y) $$
when $d \geq 3$, and
$$H(x,y) = \sum_{n=0}^\infty [\P (SRW_x (n) = y) - \P (SRW_x (n) = x) ] $$
when $d = 1$ or $2$.  It is easy to see that the sums are finite and 
that $H(x,y)$ is harmonic in $y$ except at $y=x$; it will also 
be shown that $H$ is symmetric.  [Later, we will use the above
definition of $H$ for finite graphs as well; see the appendix.]
Now let $\V0 \subseteq \VG$ be the
subspace of all functions with finite support.  Think of the vertices
of $G$ as embedded in $\V0$ by $x \mapsto \dd_x$ and the oriented
edges $\xy$ as embedded n $\V0$ by $\xy \mapsto \dd_x - \dd_y$.  
Now extend $H$ to a bilinear map on $\V0 \times \V0$.  Similarly,
think of the functions $g^u$ from Theorem~\ref{exhibit Green},
or in general the solution $g^f$ to $(I-A) g = f$ from
Corollary~\ref{general Green} as defining a bilinear form $g$ on $\V0 
\times \V0$ (or when $d \leq 2$, on part of $\V0 \times \V0$) by 
letting $g (f,\dd_x) = g^f (x)$ and extending linearly.  We then have

\begin{th} \label{g=H}
$g = H$ whenever $g$ is defined.  Consequently, $g$ is symmetric.  $\Cox$
\end{th}
 
Define the {\em transfer impedance} of two oriented edges $e$ and $f$
to be $g(e,f)$.  For any finite set $e_1 , \ldots , e_k$ of
edges, define their {\em transfer impedance matrix} $M = M(e_1 , \ldots ,
e_k)$ to be the $k$ by $k$ matrix with $M(i,j) = D^{-1} g(e_i , e_j)$,
where $D$ is the degree of the graph $G$.  Observe that the determinant
of the transfer impedance matrix is independent of the orientation
of the edges, since changing the orientation of $e_i$ has the effect
of multiplying both the $i^{th}$ row and the $i^{th}$ column of
$M$ by $-1$.  

\begin{th} \label{transfer}
Let $G$ be any periodic graph satisfying the assumptions of the 
first section.  For $e_1 , \ldots , e_k$ edges of $G$, pick an 
orientation for each edge and let $M$ denote their transfer 
impedance matrix, so $D M(i,j) = g(e_i , e_j) = 
g^{\dd_x - \dd_y} (z) - g^{\dd_x - 
\dd_y} (w)$, where $e_i = \xy$ and $e_j = \zw$.  (Here $g$ may be defined
by Corollary~\ref{general Green} or by extending $H$ linearly,
if $H$ is already known.)
If $\Tree$ is a uniform essential spanning forest for $G$, then 
$$\P (e_1 , \ldots , e_k \in \Tree) = \det (M) .$$
\end{th}

Here is an outline of why Theorem~\ref{transfer} is true.  For an 
oriented edge $e = \xy$, the function
 $(1/D) g (e, \cdot)$ gives the voltages at vertices of $G$ when 
each edge is a one-ohm resistor and one amp of current is run
from $x$ to $y$.  A straight-forward
application of Cramer's rule then produces a function $h (z)$
that is a linear combination of $g (e_i , \cdot )$ for
$i = 1 , \ldots , k$ and which computes the voltage for a unit
current across $e_k$ in the graph $\Z^d / e_1 , \ldots , e_{k-1}$.
By the equivalences between electrical networks, random walks
and uniform spanning trees \cite{Pe}, $h (x) - h(y)$ computes the
probability $\P (e_k \in \Tree \| e_1 , \ldots , e_{k-1} \in \Tree )$.
The conditional probability turns out to be $\det (M (e_1 , \ldots ,
e_k)) / \det (M (e_1 , \ldots , e_{k-1}))$ and multiplying these
together gives $\P (e_1 , \ldots , e_k \in \Tree ) = \det (M)$.

\begin{lem} \label{gives answer}
Let $G$ be a periodic graph and for edges $e_1 , \ldots , e_k$ 
forming no loop let $G' = G / e_1 , \ldots , e_{k-1}$.  Let
$e_k = \xy$ in $G$.  Let $\phi$
be a bounded function on the vertices of $G'$ harmonic 
everywhere except at $x$ and $y$, with excess
$1/D$ at $x$ and $-1/D$ at $y$.   If $\Tree$ is
the uniform random essential spanning forest of $G$ and $\Tree'$ is
the uniform random essential spanning forest of $G'$, then
$$\P (e_k \in \Tree') = \P (e_k \in \Tree \| e_1 , \ldots e_{k-1}
    \in \Tree) = \phi(x) - \phi(y) .$$
\end{lem}

\noindent{Proof:}  The first equality is standard \cite{Pe}.
For the other one, recall from \cite{Pe} that $\P (e_k
\in \Tree')$ is defined as the limit of $\P (e_k \in
\Tree_n')$ where $\Tree_n'$ is the uniform random spanning
tree on $G_n' = G_n / e_1 , \ldots , e_{k-1}$.  This probability is just
the probability that $SRW_x^{G_n'}$ first hits $y$ by moving from
$x$.  Now Corollary~\ref{converges} shows that the probability
$h(x,y,\cdot,n)$ of $SRW_\cdot^{G_n'}$ first hitting 
$y$ from $x$
converges as $n \rightarrow \infty$ to some function $h(x,y,z)$.
Since $h(x,y,z,n)$ is harmonic in $z$ except at $x$ and 
$y$, so is the limit.
(In this notation, the probability we are after is $h(x
,y,x)$.)
Also, it is easy to see that the excess of $h(x,y,\cdot,n)$ at 
$x$ is $1/D$ (use the forward equation).  Thus the excess of
$h(x,y,\cdot , n)$ at $y$ must be $-1$
and the limit satisfies
$(I-A) h(x,x,\cdot) = (\dd_{x} - 
\dd_{y})/D$.  By Corollary~\ref{unique},
there is only one such function up to an additive constant, 
so $\phi$ must equal $h(x,y,\cdot)$,
thus $\phi(x) - \phi(y) = h(x,y,x) - h(x,y,y) = h(x,y,x)$.  $\Cox$

\noindent{Proof of} Theorem~\ref{transfer}:  Proceed by induction on $k$.  
When $k = 1$, $\det M = M(1,1) = [g^{\dd_x - \dd_y} (x) - 
g^{\dd_x - \dd_y} (y)]/D$, where $e_1 = \xy$.  Denote this by 
$(1/D)g^{e_1}$.  According to Theorem~\ref{general Green},
$(1/D)g^{e_1}$ is bounded and solves $(I-A)g = [\dd_x - \dd_y]/D$,
so by the previous lemma, $M(1,1)$ calculates the probability
of $e_1 \in \Tree$.  

Now assume for induction that the theorem is true for $k-1$.  The
easy case to dispose of is when $M' \deq M(e_1 , \ldots , e_{k-1})$
has zero determinant.  Then by induction
$\P (e_1 , \ldots , e_{k-1} \in \Tree) = 0$ and it follows (e.g. from
the random walk construction of $\Tree$ in \cite{Al2,Br,Pe}) that
the edges $e_1 , \ldots , e_{k-1}$ form some loop.  Suppose without 
loss of generality that the loop is given by $e_1 , \ldots , e_r$
and that all edges are oriented forward along the loop.  Then
$\sum_{i=1}^r g^{e_i} = 0$,
hence $M$ is singular and $\det (M) = \P (e_1 , \ldots , e_k \in \Tree)
= 0$.  

In the case where $\det (M') \neq 0$, the inductive hypothesis
says that $\P (e_1 , \ldots , e_{k-1} \in \Tree) = \det (M')$
and it therefore suffices to show
\begin{equation} \label{tobeshown}
\P (e_k \in \Tree \| e_1 , \ldots , e_{k-1} \in \Tree ) = \det (M)
    / \det (M') .
\end{equation}

To show (\ref{tobeshown}) we construct the function $\phi$
of Lemma~\ref{gives answer}.  Write ${\vec {x_i y_i}}$ for $e_i$.
Electrically, what we will be doing is starting with
the function $(1/D)g^{e_k}$, which is the voltage function for
one unit of current put in at $x_k$ and taken out at $y_k$, and
adjusting it by adding a linear combination of functions
$g^{e_j}$ for $j < k$ in order to exactly cancel the 
current through each $e_j , j < k$.  This is then the voltage function
for $G/e_1 , \ldots , e_{k-1}$ when one unit of current is run 
across $\pi (e_k)$, and thus its difference across
$e_k$ computes $\P (e_k \in \Tree')$.  

To do this formally, let $\alpha_i$ for $i = 1 , \ldots , k-1$ be
real numbers for which $M_{ik} + \sum_{i=1}^{k-1} \alpha_i M_{ij} = 0$ for
$1 \leq j \leq k-1$.  These equations uniquely define the $\alpha_i$
because the columns of $M'$ are linearly independent and hence there
is a unique linear combination of them summing to $v$ where $v_i
= - M_{ik}$.  Let $N$ be the $k$ by $k$ matrix for which $N_{ij} = 
M_{ij}$ for $j < k$ and $N_{ik} = M_{ik} + \sum_{j=1}^{k-1} 
\alpha_j M_{ij}$; in other words, the first $k-1$ columns are
used to zero the nondiagonal elements of the last column of $M$
and $N$ is the resulting matrix.  Then $\det M = \det N = N_{kk}
\det (N_{ij} : i,j \leq k-1) = N_{kk} \det M'$, whence 
$N_{kk} = \det M / \det M'$.  Now define
a bounded element $J$ of $\VG$ by
$$J = D^{-1} \sum_{i=1}^k \alpha_i g^{e_i} , $$
with $\alpha_k \deq 1$.  We verify that
\begin{quote}
(i) For $j < k$, $J(x_j) - J(y_j) = 0$; \\
(ii) $J(x_k) - J(y_k) = \det (M) / \det (M')$.
\end{quote}
Indeed, $J(x_j) - J(y_j) = D^{-1} \sum_{i=1}^k \alpha_i g^{e_i} (x_j)
- g^{e_i} (y_j) = D^{-1} \sum_{i=1}^k \alpha_i DM_{ij}$ so $(i)$ 
follows from the definition of $\alpha_i$ and $(ii)$ follows from
the determination of $N_{kk}$ above.  

Now we have just shown that $J$ is constant on the pre-image of any
vertex under the contraction map $\pi$, and hence there is a
well-defined function $\phi$ on the vertices of $G$ for which $J = 
\phi \circ \pi$.  The excess of $\phi$ at a point $z$ is just the
sum over edges $\zw$ in $G$ of $\phi(z) - \phi(w)$.  This is just
the sum of $J(u) - J(v)$ over edges $\uv$ for which $\pi (\uv ) 
= \zw$, which is just the sum of excess of $J$ at $u$ over
$u \in \pi^{-1} (z)$.  The excess of $J$ at $u$, $(I-A)J(u)$, is just
$\sum_i \alpha_i (I-A) D^{-1} g^{e_i} (u)$
which is just $D^{-1} \sum_i \alpha_i (\dd_{x_i} (u) - \dd_{y_i}
(u))$.  This must be summed over $\pi^{-1} (z)$, which, for
$i < k$ contains $x_i$ if and only if it contains $y_i$.  Then
the only possible nonzero contribution to the sum is
$\alpha_k (\dd_{x_k} (u) - \dd_{y_k} (u))$ and summing this
over $u \in \pi^{-1} (z)$ gives $1$ if $z=\pi (x_k)$, $-1$ if
$z = \pi (y_k)$ and zero otherwise.  Thus $\phi$ satisfies the
conditions of Lemma~\ref{gives answer} and hence
$\P (e_k \in \Tree \| e_1 , \ldots , e_{k-1} \in \Tree ) 
= \phi(\pi (x_k)) - \phi(\pi (y_k)) = J(x_k) - J(y_k) = \det (M) / \det (M')$ 
by property (ii) above.  This finishes the induction and the proof
of the theorem.    $\Cox$

For ease of calculation, we derive a corollary to Theorem~\ref{transfer}.
\begin{cor} \label{incl-excl}
With $\Tree , M$ and $e_1 , \ldots , e_k$ as in Theorem~\ref{transfer}, pick
an integer $r$ with $0 \leq r \leq k$.  Define a $k$ by $k$ matrix
$M^{(r)}$ by $M^{(r)} (i,j) = M(i,j)$ if $i > r$ and $\dd_{ij} - M(i,j)$
if $i \leq r$.  Then $\P (e_1 , \ldots , e_r \notin \Tree, \,
e_{r+1} , \ldots , e_k \in \Tree ) = \det (M^{(r)})$.
\end{cor}

\noindent{Proof:}  The assertion for $r = 0$ is just the previous theorem; 
now assume for induction it is true for $r-1$.  By linearity of the
determinant in each row of a matrix, we have $\det M^{(r)} +
\det M^{(r-1)} = \det P$ where $P(i,j) = M^{(r)} (i,j) = M^{(r-1)} (i,j)$
if $i \neq r$ and $P(i,j) = \dd_{ij}$ if $i = r$.  Expanding by minors
along the $r^{th}$ row of $P$ gives $\det P = \det M^{(r-1)} (e_1 , \ldots ,
e_{r-1} , e_{r+1} , \ldots , e_k)$ which is equal to $\P (e_1 ,
\ldots , e_{r-1} \notin \Tree , \, e_{r+1} , \ldots , e_k \in \Tree)$
by induction.  Also by induction, $\det M^{(r-1)} = \P (e_1 ,
\ldots , e_{r-1} \notin \Tree , \, e_r , \ldots , e_k \in \Tree)$,
whence by subtraction, $\det M^{(r)} = \P (e_1 , \ldots , e_r 
\notin \Tree , \, e_{r+1} , \ldots , e_k \in \Tree)$, as desired.  $\Cox$

We end this section with a discussion of the case $d=0$, or in
other words, finite graphs.  Let $G$ be a connected
aperiodic (i.e. non-bipartite) graph on the vertices $S = 
\{ 1 , \ldots , k \}$, where as usual self-edges
have been added to make the graph $D$-regular.  The transition 
matrix $A$ has a simple eigenvalue of $1$, so it is
immediate that for any $u$ such that
$\sum_i u_i = 0$ there is a solution $g^u \in \VS$ to $(I-A) g^u = u$
and it is unique up to an additive constant.  Defining the
transfer impedance matrix by $M(\xy ,\zw ) = g^{\dd_x - \dd_y} (z)
- g^{\dd_x - \dd_y} (w)$ as before, Lemma~\ref{gives answer}
shows again that $M(1,1)$ calculates $\P (e_1 \in \Tree )$ 
and the induction is completed as before, showing that
$\P (e_1 , \ldots , e_k \in \Tree) = \det M(e_1 , \ldots , e_k)$.
Now add another self-edge to each vertex so the degrees are
all $D+1$.  If $A'$ is the new transition matrix then 
$I-A' = (1-(D+1)^{-1}) (I-A)$ and hence the solution $g^u{}'$ to
$(I-A') g^u{}' = u$ is just $(1 + D^{-1}) g^u$.  
Thus the new transfer impedance matrix is the same as the old one,
and hence the transfer impedance matrix is independent of the 
degree $D$ at which we choose to equalize the loops.  The electrical
explanation for this is that $M(e,f)$ is the induced voltage across
$f$ for a unit current with source $x$ and sink $y$, where
$e = \xy$ and every edge of $G$ is a one ohm resistor.  (To prove
this just add self-edges to regularize the degree; this
leaves $M$ and the electrical properties unchanged and they now
solve the same boundary value problem.)  The random walk interpretation
\cite{DS} is that $M(e,f)$ is the expected number 
of signed transits across $f$ of $SRW_x$ stopped when it hits $y$.

In particular, suppose $G$ is a finite, connected graph but 
not necessarily having vertices of the same degree.  We have
seen that the transfer impedances for $G$ may be unambiguously defined 
as the transfer impedances for any graph that extends $G$ to
a $D$-regular graph for some $D$ by addition of self-edges.
Write $deg_* (x)$ for the number edges incident to
$v$ that are not self-edges.  Then $deg_*$ is invariant under
degree equalization.  The relevance of $deg_*$ to transfer impedances
is that for many graphs $M(e,f)$ is approximately equal to
$\sum_{x \in e \cap f} deg_* (x)^{-1}$.  In other words, 
$M(\xy,\xy)$ is aproximately $deg_* (x)^{-1} + deg_* (y) ^{-1}$,
$M (e,f)$ is approximately $deg_* (x)^{-1}$ if $e$ and $f$ intersect at $x$,
and $M(e,f)$ is approximately zero in other cases.  

To see why this should be true, consider the condition
\begin{equation} \label{mixing}
\begin{array}{ll}
deg_*(x) \geq 2\ee^{-1} \hspace{1in} \mbox{  and  } \\[3ex]
(1 - \ee) \displaystyle{ {deg_* (x) \over deg_* (x) + deg_* (y)} 
    \leq \P (SRW_z \mbox{ hits $x$ before } y) \leq (1+\ee) {deg_* (x) 
    \over deg_* (x) + deg_* (y)} }  \end{array}
\end{equation} 
for all distinct $x,y,z$ in some subset $W$ of $G$.  This
says that a SRW from any point other than $x$ or $y$ in $W$ will
get ``lost'' with high probability before hitting $x$ or $y$
and will in fact end up choosing which to hit first in proportion
to the stationary measures at the two vertices.  There are many
families of graphs $G_k$ for which~(\ref{mixing}) holds with $W=G$
for a sequence $\ee_k$ converging to zero as $k \rightarrow \infty$;
examples include the complete graph on $k$ vertices and
the $k$-cube.

Assume condition~(\ref{mixing}) for some $\ee > 0$.  For some vertex
$x$, all of whose neighbors are in $W$, enumerate its neighbors
$y_1 , \ldots , y_{deg_* (x)}$
(allow repeated neighbors if there are parallel edges).  
Let $e_i = \xy_i$ for $i \leq deg_* (x)$.  Use the interpretation
of $M(e_i , e_j)$ as the expected number of signed transits across
$e_j$ for a random walk started at $x$ and stopped at $y_i$.  Reversibility
implies that the expected number of signed transits is zero over 
all times before the last visit to $x$, so 
conditioning on the first step of $SRW_x$ and using the
``craps'' principle shows that for fixed $i$ and varying $j$,
$M(e_i , e_j)$ is proportional to $\P (SRW_{y_j} \mbox{ hits $y_i$ before }
x)$.  Applying condition~(\ref{mixing}) together with the fact that
$\P (SRW_{y_i} \mbox{ hits } y_i \mbox{ before } x) = 1$ shows
that, when $j \neq i$,
\begin{eqnarray*}
&& {(1 - \ee) deg_* (y_i) / (deg_* (x) + deg_* (y_i)) \over
1 + (deg_* (x) - 1) (1 + \ee) deg_* (y_i) / (deg_* (x) + deg_* (y_i)) } \\[2ex]
\leq M(e_i,e_j) & \leq &
   {(1 + \ee) deg_* (y_i) / (deg_* (x) + deg_* (y_i)) \over
   1 + (deg_* (x) - 1) (1 - \ee) deg_* (y_i) / (deg_* (x) + deg_* (y_i)) } .
\end{eqnarray*}
Since $deg_* (x) , deg_*(y) \geq 2\ee^{-1}$, $deg_*(x) deg_*(y) / (
deg_*(x) + deg_*(y)) \geq \ee^{-1}$ and the addition of $1$ in the
denominator of the first term in the above inequality loses no
more than a factor of $1 -\ee$, and we may rewrite the inequalities as
$$ (1-\ee)^3 deg_*(x)^{-1} \leq M(e_i, e_j) \leq (1-\ee)^{-2} deg_* (x)^{-1} .$$
Similarly, 
$$ (1-\ee)^3 (deg_*(x)^{-1} + deg_* (y_i)^{-1}) \leq M(e_i, e_i) \leq 
   (1-\ee)^{-2} (deg_* (x)^{-1} + deg_* (y_i)^{-1}) .$$
Finally, for $e = \xy$ and $f =\zw$ such that $x,y,z,w$ and all
the neighbors of $x$ are in $W$, we have $M(e,f) = M(e,e) [
\P (SRW_z \mbox{ hits $x$ before }y) - \P (SRW_w \mbox{ hits $x$ before }y)]
\leq 2\ee (1 - \ee)^{-2} (deg_* (x)^{-1} + deg_* (y)^{-1})$.  
This follows from the electrical interpretation of $M(e,f)$, since a
unit current flow puts a voltage difference of $M(e,e)$ across $e$, 
after which the voltages elsewhere are $M(e,e)$ times the probability
from there of $SRW$ hitting $x$ before $y$.  Thus
the mixing condition~(\ref{mixing}) does indeed imply that
\begin{equation} \label{M is right}
M(e,f) = \sum_{x \in e \cap f} deg_* (x)^{-1} + O(\ee ) \sum_{x \in e}
    deg_* (x)^{-1} 
\end{equation}
for $x,y,z,w \in W$ with all neighbors of $x$ inside $W$.

\section{An Example and a High Dimensional Limit}

The case where $G$ is the nearest neighbor graph for $\Z^2$ is
special because the Green's function can be explicitly evaluated
as a polynomial in $\pi^{-1}$.  Following \cite{Sp}, we have
that the Green's function is given by
$$ H(0,x) = (2 \pi)^{-2} \int {1 - \cos (x \cdot \alpha) \over 
  1 - (1/2) \cos (\alpha_1) - (1/2) \cos (\alpha_2)} \, d\alpha  $$
and for $x = (n,n)$ a change of variables from $(\alpha_1 , \alpha_2)$
to $(\alpha_1 + \alpha_2 , \alpha_1 - \alpha_2)$ yields
$$H((0,0),(n,n)) =  4 \pi^{-1} \left [ 1 + {1 \over 3} + \cdots + 
{1 \over 2n-1} \right ].$$
These values, along with the symmetries of the lattice and the fact 
that $H$ is harmonic, allow $H$ to be determined recursively, the 
first few values being 

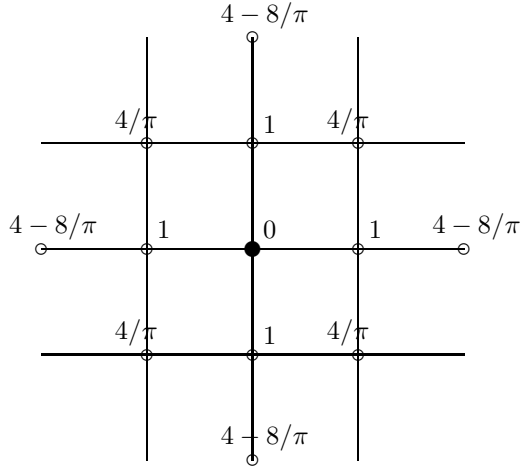
\setlength{\unitlength}{2pt}
\begin{picture}(80,120)(-80,-20)
\put(40,40) {\circle* {3}}
\put(42,42) {0}
\put(40,60) {\circle {2}}
\put(42,62) {1}
\put(40,80) {\circle {2}}
\put(34,83) {$4-8/\pi$}
\put(40,20) {\circle {2}}
\put(42,22) {1}
\put(40,0) {\circle {2}}
\put(34,3) {$4-8/\pi$}
\put(60,40) {\circle {2}}
\put(62,42) {1}
\put(80,40) {\circle {2}}
\put(74,43) {$4-8/\pi$}
\put(20,40) {\circle {2}}
\put(22,42) {1}
\put(0,40) {\circle {2}}
\put(-6,43) {$4-8/\pi$}
\put(60,60) {\circle {2}}
\put(54,63) {$4/\pi$}
\put(60,20) {\circle {2}}
\put(54,23) {$4/\pi$}
\put(20,20) {\circle {2}}
\put(14,23) {$4/\pi$}
\put(20,60) {\circle {2}}
\put(14,63) {$4/\pi$}
\put(20,0) {\line(0,1){80}}
\put(40,0) {\line(0,1){80}}
\put(60,0) {\line(0,1){80}}
\put(0,20) {\line(1,0){80}}
\put(0,40) {\line(1,0){80}}
\put(0,60) {\line(1,0){80}}
\end{picture}

\noindent{Let} $w_1, w_2, w_3, w_4$ denote respectively the edges connecting
the origin to $(1,0),(0,1),(-1,0),(0,-1)$.  The above values for
$H$ then yield the following circulant for $M$:
$$M(w_1 , \ldots , w_4) = \left [ \begin{array}{cccc} 
    1/2 & 1/2 - \pi^{-1} & 2\pi^{-1} - 1/2 & 1/2 - \pi^{-1} \\
    1/2 - \pi^{-1} & 1/2 & 1/2 - \pi^{-1} & 2\pi^{-1} - 1/2 \\
    2\pi^{-1} - 1/2 & 1/2 - \pi^{-1} & 1/2 & 1/2 - \pi^{-1} \\
    1/2 - \pi^{-1} & 2\pi^{-1} - 1/2 & 1/2 - \pi^{-1} & 1/2 
  \end{array} \right ] . $$

Theorem~\ref{transfer} and Corollary~\ref{incl-excl} assert for example,
that $\P (w_1 , w_2 , w_3 , w_4 \in \Tree)$ and $\P (w_1 , w_2 , 
w_3 \notin \Tree , \, w_4 \in \Tree)$ are given respectively by 
$\det M$ and $\det M^{(3)}$ respectively, where
$$ M^{(3)} = 
    \left [ \begin{array}{cccc} 
    1/2 & -1/2 + \pi^{-1} & -2\pi^{-1} + 1/2 & -1/2 + \pi^{-1} \\
    -1/2 + \pi^{-1} & 1/2 & -1/2 + \pi^{-1} & -2\pi^{-1} + 1/2 \\
    -2\pi^{-1} + 1/2 & -1/2 + \pi^{-1} & 1/2 & -1/2 + \pi^{-1} \\
    1/2 - \pi^{-1} & 2\pi^{-1} - 1/2 & 1/2 - \pi^{-1} & 1/2 
\end{array} \right ] . $$
The determinants of $M$ and $M^{(3)}$ are respectively
$(4 \pi^{-1} - 1) (2 \pi^{-1} - 1)^2$ and $2 \pi^{-2} - 4\pi^{-3}$,
so the probability of all four edges incident to $0$ being in
$\Tree$ is $(4 \pi^{-1} - 1) (2 \pi^{-1} - 1)^2 \approx .0361$, 
while the probability that the origin is a leaf of $\Tree$ (i.e. has 
degree one in $\Tree$) is $8 \pi^{-2} - 16\pi^{-3} \approx .2945$.

The remainder of this section corrects, proves and generalizes some
conjectures of Aldous about spanning trees for graphs as the graphs
tend to infinity {\em locally}, in the sense that the minimum number of
neighbors of a vertex all grow without bound.  We first quote Conjecture 11 
from \cite{Al2}.  To do this, let $G_k$ denote a sequence of finite 
graphs, each with a distinguished vertex $v_k$.  Let $r_k = deg_* (v_k)$,
let $A_k$ be the set of neighbors
of $v_k$ in $G_k$ and for $w \in A_k$, let $\psi_k (w) = 
\P (SRW_w^{G_k} \mbox{ hits $A_k \setminus \{ w \}$ before } v_k)$.  
Aldous then conjectures the following \cite[Conjecture 11]{Al2}.
\begin{quote}
{\em Let $1 + D_k$ denote the random degree of $v_k$ in the uniform
random spanning tree on $G_k$.  Suppose that
$$ r_k \rightarrow \infty \hfill ; \hfill \sup_w |r_k (1 - 
    \psi_k (w)) - 1| \rightarrow 0 .$$
Then $D_k$ converges in distribution to a Poisson with mean 1.  }
\end{quote}

Here is a counterexample to the conjecture.  Let the vertices of
$G_k$ other than $v_k$ be $\{ x_i , y_i , z_{ij} : 1 \leq i \leq k ; 
1 \leq j \leq 4k \}$, with edges connecting $v_k$ to each $x_i$
and $y_i$ and for every $i$, an edge connecting $x_i$ to each $z_{ij}$ 
and an edge connecting $y_i$ to each $z_{ij}$.  Then $r_k = 2k$.
By symmetry, we have $\psi_k (w) = \psi_k (x_1)$ for any $w \in A_k$.
This is equal to $\P (SRW_{x_1}
\mbox{ hits $y_1$ before } v_k) = 2k/(1+2k)$, hence 
$r_k (1 - \psi_k (w)) = 2k/(2k+1) \rightarrow 1$ and the hypothesis
in the conjecture is satisfied.  But for each $i$, any spanning
tree contains either an edge connecting $v_k$ to $x_i$ or an
edge connecting $v_k$ to $y_i$, so the degree of $v_k$ is at least $k$.

Evidently, a condition different from $r_k (1 - \psi_k (w))$ 
converging uniformly to 1 is required for the conjecture to
be true.  Aldous' condition is trying to capture two aspects of the
graph: some sort of mixing (SRW from any neighbor of $v_k$ returns to 
$v_k$ before $A$ with the same probability) and the correct
total probabilities (these probabilities should all be about
$r_k^{-1}$ so they can sum to 1).  The mixing part of the
condition as stated in the conjecture is too weak, as illustrated
by the counterexample, and needs to be replaced by a condition
that equalizes the individual return probabilities of $SRW_w$
hitting $v_k$ before $z$ for any neighbors $w,z$ of $v_k$.  The 
most natural such condition from our viewpoint is~(\ref{mixing})
with $x = v_k$ and $\ee = \ee_k$ for some $\{ \ee_k \}$ going
to zero.  On the other hand, there is no need to require all the
neighbors of $v_k$ to have the same degree.  Once a sufficient
amount of independence has been achieved,
it suffices for the expected number of such edges in the tree to be 
converging to a constant, $1+\lambda$.  The conjecture may thus be 
revised to yield the following theorem.
\begin{th} \label{trueconj11}
Let $G_k$ be a sequence of finite graphs.  Let $\ee_k \rightarrow 0$ be
a sequence of positive numbers and let $v_k, A_k, D_k$ and $r_k$ be as 
above.  Assume that (\ref{M is right}) holds with $W = A_k$ and 
$\ee = \ee_k$; this is implied for example by~(\ref{mixing}).  
If in addition, $\sum_{w \in A_k}
1 / deg_* (w) \rightarrow 1+\lambda$ as $k \rightarrow \infty$
for some positive $\lambda$, then $D_k$
converges to a Poisson with mean $\lambda$.
\end{th}

\noindent{Remark:}  Since the theorem gives local behavior at $v_k$,
it can easily be extended to the case where $G_k$ are infinite
graphs on which there is a Harnack principle: simply take $G_k'$ 
to be a large enough finite piece of $G_k$ so that the hypotheses
of the theorem are true (this is possible by the Harnack principle).

The following lemma will be necessary when calculating determinants
of transfer impedance matrices.
\begin{lem} \label{detcalc}
Let $a_1 , \ldots , a_{k+1}$ be positive real constants and let $e_1 , \ldots
e_k$ be the edges of a spanning tree whose vertices are
$\{ 1 , \ldots , k+1 \}$.  Define a $k$ by $k$ matrix $M$ by
letting $M(i,j)$ be $a_r + a_s$ if $i = j$ and $e_i$ connects
$r$ to $s$; $a_r$ if $e_i$ and $e_j$ are distinct edges meeting at 
$r$; and zero otherwise.  Then $\det M = (\prod a_i) (\sum a_i^{-1})$. 
\end{lem}

\noindent{Proof:}  If $k=1$ then $M = (a_1 + a_2)$ and the lemma
is clearly true.  now assume for induction that the lemma is true
for $k-1$.  Assume by renumbering if necessary that the vertex $1$ is
a leaf of the tree, there being a single edge $e_1$ connecting
$1$ to $2$.  Also assume that $2$ is connected to $1,3,4, \ldots , r$
by $e_1 , \ldots , e_{r-1}$ respectively, 
for some $r \leq k+1$.  Row reduce $M$ by subtracting $a_2 / (a_1 +
a_2)$ times the first row from rows $2, \ldots, r$.  Since edges $2,
\ldots , r$ were the only edges incident to $e_1$, this clears the first
column.  The remaining entries of $M$ are unchanged except that
all appearances of $a_2$ get replaced by $a_1 a_2 / (a_1 + a_2)$.
Expanding along the new first column gives $\det M = (a_1 + a_2)$ times
the determinant of the $k-1$ by $k-1$ matrix gotten by taking all but
the first row and column of $M$ and replacing $a_2$ by $a_1 a_2 /
(a_1 + a_2)$.  By induction, the latter determinant is $(a_1 / 
(a_1 + a_2)) (\prod_{i \geq 2} a_i) ((a_1 + a_2) / a_1 a_2 +
\sum_{i \geq 3} a_i^{-1})$, so $\det M = (\prod a_i) ((a_1 + a_2) / a_1
a_2 + \sum_{i \geq 3} a_i^{-1}) = \prod a_i \sum a_i^{-1}$ as desired.   $\Cox$ 

\noindent{Proof} of Theorem~\ref{trueconj11}:  Let $X$ be a
random variable for which $X-1$ is Poisson with mean $\lambda$. 
Set
$$\phi (z) = \sum \P (X = n) z^n = z e^{\lambda (z-1)} .$$
Then the $s^{th}$ factorial moment $\E (X)_s \deq \E X(X-1) \cdots (X-s+1)$ 
of $X$ is the $s^{th}$ derivative of $\phi$ at 1, which 
is equal to $\lambda^s+s\lambda^{s-1}$.
The factorial moments determine this distribution uniquely,
since its moment generating function exists in a neighborhood of
zero, from which it follows that convergence of the factorial
moments of a sequence of random variables to the factorial
moments of $X$ implies convergence in distribution to $X$
(see for example Theorem 3.10 of~\cite{Du}).
It suffices therefore to prove that $\E (1+D_k)_s \rightarrow 
\lambda^s+s\lambda^{s-1}$ as $k \rightarrow \infty$ for each $s$.  

Write $\E (1+D_k)_s$ as $\sum \P (e_1, \ldots , e_s \in \Tree )$
where the sum is over all ordered collections $(e_1 , \ldots , e_s)$
of $s$ distinct edges incident to $v_k$.  Fixing such a
collection, we have $\P (e_1 , \ldots , e_s \in \Tree) = 
\det M(e_1 , \ldots , e_s)$.  From~(\ref{M is right}) we have
that $M(i,j) = (1 + O(\ee_k)) (deg_* (x)^{-1} + \dd_{ij} deg_* 
(y_i)^{-1})$ where $e_i = \xy_i$.  With $s$ staying fixed, 
this and Lemma~\ref{detcalc} give
\begin{equation} \label{sum}
\det M = (1+o(1)) (\prod_{i=1}^s
deg_* (w_i)^{-1}) (1 + r_k^{-1} \sum_{i=1}^s deg_*  (w_i)) .
\end{equation}
Summing over all ordered collections gives 
$$ \det M = (1+o(1)) [\sum_{(e_1 , \ldots , e_s)} \prod_{i=1}^s 
   deg_* (w_i)^{-1} + \sum_{(e_1 , \ldots , e_{s-1})} s (r_k - s) 
   r_k^{-1} \prod_{i=1}^{s-1} deg_* (w_i)^{-1}] $$
since each ordered collection of size $s-1$ appears $s (r_k-s)$ times
in the second term of~(\ref{sum}) and each ordered collection of
size $s$ appears once in the first term of~(\ref{sum}).  As $s$
remains fixed with the minimum degree among $x$ and its neighbors
converging to infinity, the above expression for $\det M$ converges
to $(\sum_i deg_* (w_i)^{-1})^s + s(\sum_i deg_*  (w_i)^{-1})^{s-1} \rightarrow
\lambda^s + s\lambda^{s-1}$, proving the convergence of factorial moments 
and the theorem.  $\Cox$

Suppose now that we are interested not only in the degree of $v_k$
but in the local structure of the essential spanning forest near
$v_k$.  For any locally finite rooted tree $\Tree$, let $\Tree 
\wedge r$ denote the random finite subtree of $\Tree$
consisting of vertices connected to the root by paths of length at
most $r$ in $\Tree$.  We say that a sequence of tree-valued random
variables converges in distribution (written $\Tree_k
\dconv \Tree$) if $\Tree_k \wedge r \dconv \Tree \wedge r$ for
every $r$, where the latter is defined to hold when $\P (\Tree_k \wedge
r = t) \rightarrow \P(\Tree \wedge r = t)$ for every $t$ of height at
most $r$.  Under suitable conditions on the graphs $G_k$, it will
turn out that the component $\Tree_k$ of the uniform essential spanning
forest on $G_k$ rooted at $v_k$ will converge in distribution 
to a particular tree $\pois$ which we now define.
Let $\pois$ be a singly infinite path,
$x_0 , x_1 , \ldots$, to which has been added at each
$x_i$ the tree of an independent Poisson (1) branching process
(which is critical hence finite with probability one).  Another
way of describing $\Tree$ is as the tree of a Poisson (1) branching
process rooted at $x_0$ and conditioned to survive forever. 
Aldous has conjectured (personal communication, though the
conjecture is implicit in \cite{Al3}) that 
$\Tree_k$ converges in distribution to $\pois$ whenever $G_k$
grows locally in a sufficiently regular manner.  In the
terminology of \cite{Al1}, $\Tree_k$ should converge to a sin-tree
with the fringe distribution of a Poisson (1) branching process.
This is known in the special case where $G_k$ = $K_k$, the complete
graph on $k$ vertices \cite{Al3,Gr}.

We are now in a position to prove this.  A consequence is that
the probability of there existing two disjoint paths of length $L_k$
in $\Tree_k$ from $v_k$ goes to zero as $k \rightarrow \infty$,
provided that $L_k \rightarrow \infty$.
An question left open in \cite{Pe} is whether the components
of the uniform essential spanning forest have one or two 
ends (the possibility of more than two is ruled out by an
argument in \cite{BK}).  We believe the answer to be that
all components have one end, and convergence to zero of the
probability of there being two disjoint infinite paths from $v_k$ 
can be viewed as a heuristic argument in favor of all
components having one end.

\begin{th} \label{poisson 1}
Let $G_k$ be a sequence of finite graphs with distinguished vertices
$v_k$.  Assume, by renumbering if necessary, that 
$$(1+o(1))\max_{x \in G} deg_* (x) = k = (1+o(1)) \min_{x \in G} deg_* (x)$$
as $k \rightarrow \infty$.  Also assume the following version 
of~(\ref{M is right}) uniformly in edges $e,f \in G$:
$$ M(e,f) = k^{-1} (|e \cap f| + o(1)) . $$
Then $\Tree_k \dconv \pois$.
\end{th}

\noindent{Remark:}  For $r = 1$, $\Tree \wedge r$ is a star centered
at $v_k$ with $1 + D_k$ edges, while $\pois \wedge r$ is a star
centered at $x_0$ with $1 + X$ edges, where X is Poisson with
mean 1.  Thus the case $r=1$ is essentially the previous theorem.
Notice also that the usual families of graphs $G_k$ (e.g. complete
graph on $k$ vertices, $k$-cube, $k/2$-dimensional torus of arbitrary
length) all satisfy the hypotheses of the theorem.

\noindent{Proof:}  For a finite rooted tree $t$ and finite
rooted graph $u$, say
that a map $f$ from the vertices of $t$ to the vertices of $u$
is a tree-map if $f$ is injective, maps the root of $t$ to the root
of $u$, and $f(x) \sim f(y)$ for each $x \sim y$.  Let
$N(u;t)$ denote the number of distinct tree maps from $t$ to $u$.
For example, if $t$ and $u$ are stars of respective sizes $s$ and $r$ 
about their roots, then $N(u;t) = (r)_s$.  The proof of this theorem
generalizes the proof of the preceding theorem, in the sense that
$\E N(\Tree_k ; t)$ is a sort of generalized $t^{th}$ moment of $\Tree_k$. 
In the appendix, the usual tightness criteria for convergence of
probability measures are extended in an obvious way to tree-valued 
random variables, showing in particular (Theorem~\ref{conv criterion})
that if $\E N(\Tree_k ; t)
\rightarrow \E N(\Tree ; t)$ for each $t$ and the values of
$\E N(\Tree ; t)$ uniquely determine the distribution of $\Tree$
then $\Tree_k \dconv \Tree$.  Also proved there is the somewhat less
trivial fact (Theorem~\ref{tree moments})
that the values of $\E N(\Tree ; t)$ uniquely determine the
distribution of $\Tree$ under the growth condition: $\E N(\Tree ; t)
\leq e^{c |t|}$ for some $c$.  (A sharper growth condition such as
an analogue to Carleman's condition could be obtained but is not 
needed here.)  What remains then, is to show that $\E N(\Tree_k ; t)
\rightarrow \E N(\pois ; t)$ for each finite $t$ and to verify the
growth condition on $\E N (\pois ; t)$.

\noindent{Begin} by establishing
\begin{equation} \label{NXt}
\mbox{For any finite rooted tree } t, \E N(U;t) = 1,
\end{equation}
where $U$ is the tree of a Poisson (1) branching process rooted at
some vertex $y_0$.  Let $z_0$ be the root of $t$.  Use induction on the
height of $t$.  If $t$ is just $z_0$, then $N(u;t) = 1$
for any $u$, so the equation is trivially true.  Now suppose $z_0$
has $s$ descendants $z_1 , \ldots , z_s$ for some $s > 0$ and assume 
for induction that~(\ref{NXt}) holds for each of the subtrees $t_i, 
1 \leq i \leq s$ rooted at $z_i$.  Let $r \geq 0$ be the random
number of descendants $y_1 , \ldots , y_r$ of $y_0$.  Conditional 
upon $r$, each of the subtrees $u_i$ rooted at $y_i$
is the tree of an independent Poisson (1) branching process.  
Now any tree-map from $t$ to $U$ maps each $t_i$ to a distinct $u_j$.
Thus $N(U;t) = \sum \prod_{j=1}^s N(u_{k_j};t_j)$ where the sum
is over all ordered sequences of distinct $k_1 , \ldots , k_s$ chosen
from among $1 , \ldots , r$.  Then by independence conditional 
on $r$, 
$$\E N(U;t) = \sum_r (e^{-1} /r!) \sum \prod_{j=1}^s \E N(u_{k_j};t_j) . $$
By the induction hypothesis each expectation is one, so the sum is
just the number of ways of choosing the $k_j$'s.  Thus
$\E N(U;t) = \sum_r (e^{-1} /r!) (r)_s = \E (X)_s$ where
$X$ is a Poisson of mean 1.  This is equal to 1, establishing~(\ref{NXt}).

Now we compute $\E N(\pois ;t)$.  Recall that $\pois$ is a path
$x_i : i \geq 0$ to which has been added an independent Poisson (1)
branching process, say $U_i$ at each $x_i$.  If $f$ is a tree-map from
$t$ into $\pois$, there is a greatest $i$ for which $x_i$ is in
the range of $f$.  Let $w(f)$ denote the vertex of $t$ that maps
to $x_i$ for this greatest $i$.  For each vertex $z$ of $t$, we will
show that the expected number of tree-maps $f: t \rightarrow \pois$
for which $w(f) = z$ is one.  Indeed, if $x_0 = z_0 , z_1 , \ldots ,
z_k = z$ is the path from the root of $t$ to $z$, then the tree-maps
$f$ for which $w(f) = z$ are in one to one correspondence with
the collections of maps $f_0 , \ldots , f_k$ where $f_i$ maps
the subtree $t_i$ of $t$ rooted at $z_i$ to the subtree $\pois (i)$
 of $\pois$ rooted at $X_i$.  Thus the number of $f$ for which 
$w(f) = z$ is $\prod_{i=0}^k N(\pois (i); t_i)$.  But each $\pois
(i)$ is an independent Poisson (1) branching process, so by
equation~(\ref{NXt}), the product is one.  Finally, summing over
the vertices $z$ of $t$ gives that $\E N(\pois ; t)$ is equal to
$|t|$, the number of vertices of $t$.  This verifies the growth 
condition on $\E N(\pois ; t)$ with miles to spare!

To calculate $\E N(\Tree_k ; t)$, observe that any tree-map 
$f : t \rightarrow \Tree$ is also a tree-map from $t$ to $G_k$ 
where $G_k$ is considered to be rooted at $v_k$.  If $u(f)$ denotes
the image of $f$ as a subtree of $G_k$, then we may write
\begin{equation} \label{each prob}
\E N(\Tree_k ; t) = \sum \P (u(f) \subseteq \Tree_k )
\end{equation}
where the sum is over all tree-maps $f : t \rightarrow G_k$.
Fix $f$.  Then $\P (u(f) \subseteq \Tree_k) = \det M$ where
$M$ is the transfer impedance matrix for $u(f)$ as a subgraph of $G$.
By hypothesis $M(e,e') = k^{-1} (|e \cap e'| + o(1))$.  Since $f$ 
is injective, $f^{-1}$ is defined on $u(f)$ so $|e \cap e'| = |f^{-1} 
(e) \cap f^{-1} (e')|$, thus for any $f$ the transfer impedance 
matrix for $u(f)$ in $G_k$ may be written as $k^{-1} P$, where
$P$ is a matrix indexed by the edges of $t$ for which $P(e,e') = 
(|e \cap e'| +o(1))$ as $k \rightarrow \infty$.  Since the size
of $P$ remains fixed as $k \rightarrow \infty$, Lemma~\ref{detcalc}
with $a_i = 1$ for all $i$ gives that $\det P$ is equal to the
$|t| (1+o(1))$, and thus $\det M = k^{1 - |t|} |t| (1+o(1))$.
Then the probabilities in equation~(\ref{each prob}) are all
equal and the identity becomes 
$$\E N(\Tree_k ; t) = N(G_k , t) k^{1 - |t|} |t| (1+o(1)) .$$
But $N(G_k , t) = k^{|t|-1} (1+o(1))$.  An easy way to see 
this is to consider building a tree-map $f : t \rightarrow G_k$
starting at the root and working outwards, not worrying
about injectivity.  If $f$ is defined on
$z$ then by hypothesis there are $k (1+o(1))$ neighbors of
$f(z)$ and each descendant of $z$ may be mapped to any neighbor
of $f(z)$.  The fraction of all maps built this way that are 
injective goes to one as $t$ remains fixed and $k \rightarrow \infty$,
so the total number of maps is $k^{|t|-1} (1+o(1))$.  Now~(\ref{each
prob}) becomes $\E N(\Tree_k ; t) = |t| (1+o(1))$, hence 
$\E N(\Tree_k ; t) \rightarrow \E N(\pois ; t)$ as claimed.
Since the growth condition on $\E N(\pois ; t)$ has been verified,
this shows $\Tree_k \dconv \pois$.   $\Cox$

\section{Entropy}

In this section we consider the entropy of the essential spanning forest 
process.

The set of essential spanning forests is a closed shift-invariant subset 
of
$\{0,1\}^{E(G)}$ where $E(G)$ is the edge set of $G$.
The topological entropy (per vertex) of the essential
spanning forest is defined to be 
$$H_{top}  = \lim_{n\rightarrow\infty}{1\over |B_n|} \log (N_{B_n})$$
where $B_n$ is an increasing sequence of rectangular boxes (i.e. of the 
form
$C\times S$ where $C$ is a rectangular box in $\Z ^d$ together with the
induced edges) and $N_B$ is the number of essential spanning forests of 
the induced
graph $B$ where a forest is essential if 
every component of the graph is required to touch the boundary
of $B$.  We may also consider boxes with boundary conditions, meaning
a box $B$ together with an equivalence relation $\equiv$ on the vertices
of $B$ that neighbor $B^c$.  An essential spanning forest on a box
with boundary conditions is one that becomes a tree under the
contraction map consisting $B \mapsto B/\equiv$ (think of the boundary
conditions as telling which vertices are connected by unseen edges in
$B^c$).

Notice that if $B$ is the union of 
two boxes $C$ and $D$ then any essential spanning forest of 
$B$ restricts to essential spanning forests in both $C$ and $D$.
This means that $N_B \leq N_C N_D$ so that 
$\log (N_B)$ is subadditive and the entropy is 
independent of the sequence of boxes chosen.
There is a variational principle for the topological entropy of the essential
spanning forest.
$$H_{top} = \sup\{H(\mu) \| \mu 
\mbox{ is an invariant probability measure  }\}$$
where $H(\mu)$ is the Kolmogorov-Sinai entropy of $\mu$ per vertex with 
respect to the group of translations by $\Z ^d$.
See \cite{Mi} for a short proof of this fact.

Now let $\mu_G$ be the probability measure of uniform essential spanning 
forests on 
$G$.
Using arguments described in \cite{Pe} it is seen that if $B_n$ is any 
increasing
sequence of rectangular boxes with arbitrary boundary conditions and if
$\mu_n$ is the measure that gives equal weight to each forest in $B_n$ 
in which each component of the forest meets the boundary of $B_n$ the 
$\mu_n$
converges weakly to the translation invariant probability measure $\mu_G$,
moreover this convergence is uniform in the boundary conditions.
By uniform convergence we mean the following.
Suppose that we are given a cylinder set and an $\epsilon \geq 0$.
Then there is a box $B$, containing the cylinder,
so that for any box $C$ containing $B$ and any boundary conditions
on $C$ we have that the uniform probability measure on essential spanning 
forests 
of $C$ takes a value on the cylinder set that agrees with that of $\mu_G$ 
to within $\epsilon$.

The principle content of this section is the following theorem.
\begin{th} \label{entropy}
(a)  The measure $\mu_G$ of the uniform essential spanning forest process 
is the unique translation-invariant measure on the set of essential spanning 
forests whose Komogorov-Sinai entropy is $H_{top}$.\\
(b)$$H_{top} = \frac{1}{k} \int_{\VT} 
\log (D^k \chi ( {Q(\alpha)}(1) ))d\alpha$$
where $\chi ( {Q(\alpha)} )$ is the characteristic polynomial of $Q(\alpha)$ 
and
the integral is over the d-torus with respect to Haar measure.
\end{th}

Before proving the theorem we give some examples of $(b)$ in which the
entropy can essentially be read off.  
The case when $S = \{1\}$, i.e. $k=1$ is especially easy to analyze because 
$Q(\alpha)$ is a $1 \times 1$ matrix.  In these cases the entropies
are the same as entropies calculated by Lind et al \cite{LS} of some 
seemingly unrelated dynamical systems that can be represented as
Bernoulli shifts on certain subgroups of 
$(\Z^d)^{(\R / \Z)}$ defined by peroidic linear relations.  We are at a
loss to explain this apparent coincidence.

Suppose the origin is connected 
to $M$ pairs of opposite vertices in $G = \Z^d$.  Say the number of
self-edges per vertex is $l$, though clearly this must drop out of
the calculation.  Suppose we denote representatives 
of these pairs by $\{ x_m \, : \, 1 \leq m \leq M \} $.  Then $D = 2M + l$ 
and
$$Q(\alpha) = D^{-1} (l +  \sum _{m=1}^{M} 2 \cos( 2 \pi \alpha\cdot 
x_m )) . $$
Then $\chi (Q(\alpha))(1) = 1 - Q(\alpha) = D^{-1} (\sum_{m=1}^M 2-2\cos
(2\pi \alpha \cdot x_m))$ and the entropy is  
$$\int_{T^d} \log (2M - \sum _{m=1}^{M} 2 \cos( 2 \pi \alpha\cdot x_m 
)) d\alpha . $$
One such example is $\Z^2$ itself with nearest neighbor edge relation.  
The entropy is 
$$\int_{0}^{1} \int_{0}^{1} \log (4 - 2\cos (2\pi \alpha_{1}) - 2\cos 
(2\pi \alpha_{2} )) d\alpha_{1} d\alpha_{2}\approx 1.166 . $$
Another is the triangle lattice.  This has a representation as the nearest 
neighbor lattice on $\Z^2$ with added edges placed in the ``southwest 
- northeast'' diagonal of each square.  \\
\setlength{\unitlength}{2pt}
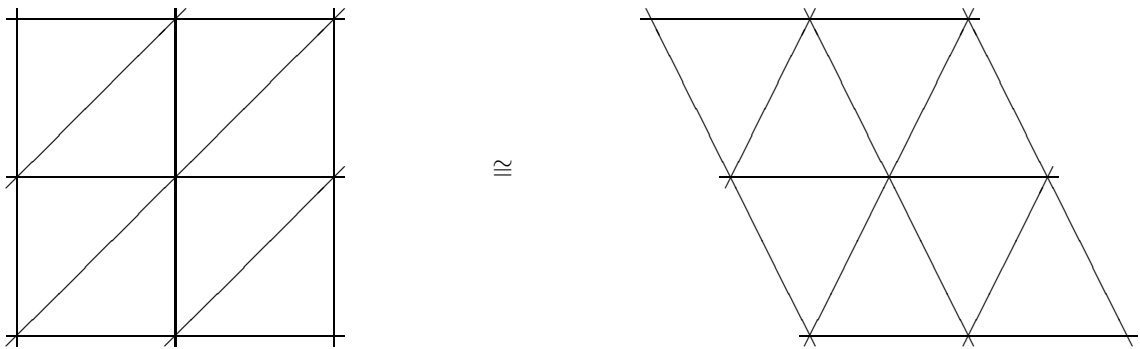
\begin{picture}(210,100)(15,0)
\put(20,18) {\line(0,1){64}}
\put(50,18) {\line(0,1){64}}
\put(80,18) {\line(0,1){64}}
\put(18,20) {\line(1,0){64}}
\put(18,50) {\line(1,0){64}}
\put(18,80) {\line(1,0){64}}
\put(18,18) {\line(1,1){64}}
\put(18,48) {\line(1,1){34}}
\put(48,18) {\line(1,1){34}}
\put(110,50){$\cong$}
\put(138,80) {\line(1,0){64}}
\put(153,50) {\line(1,0){64}}
\put(168,20) {\line(1,0){64}}
\put(171,18) {\line(-1,2){32}}
\put(201,18) {\line(-1,2){32}}
\put(231,18) {\line(-1,2){32}}
\put(169,18) {\line(1,2){32}}
\put(199,18) {\line(1,2){17}}
\put(154,48) {\line(1,2){17}}
\end{picture} \\
In this case, ignoring self-edges, $Q(\alpha) = \frac{1}{6} 
(2\cos (2\pi \alpha_{1}) + 2\cos (2\pi \alpha_{2}) + 2\cos (2\pi (\alpha_{1} 
+ \alpha_{2})))$ and the entropy is 
$$\int_{0}^{1} \int_{0}^{1} \log (6 - 2\cos (2\pi \alpha_{1}) - 2\cos 
(2\pi \alpha_{2} ) -  2\cos (2\pi (\alpha_{1} + \alpha_{2}))) d\alpha_{1} 
d\alpha_{2} \approx 1.61 . $$
These are the same as the entropies given in \cite{LS} for Haar measure
on the subgroups of $(\R / \Z)^{\Z^2}$ consisting respectively of those
configurations $\phi$ for which $4\phi (x) - \phi(x+(0,1)) - \phi(x+(0,-1)) -
\phi(x+(1,0)) - \phi(x+(-1,0)) = 0$ and those configurations $\phi$ for
which $6\phi (x) - \phi(x+(0,1)) - \phi(x+(0,-1)) - \phi(x+(1,0)) - 
\phi(x+(-1,0)) - \phi(x+(1,1)) -\phi(x+(-1,-1)) = 0$. 

The proof of Theorem~\ref{entropy} uses the following lemma on the
stability of entropy under changes in a small percentage of the 
output of the process.

\begin{lem} \label{entropy2}
Let $(\Omega , \mu )$ be a Lebesgue probability space.  
Suppose that ${X}  = (X_1 , X_2 , \ldots , X_ N )$  
and ${Y} =  (Y_1 , Y_2 , \ldots , 
Y_ N )$ are binary random variables such that for all $\omega \in \Omega
,\hspace{3pt} \# \{ i | X_i \neq   Y_i \}  \leq K.$  Then

$$ |  \frac{1}{N}  H({X}) - \frac{1}{N} H({Y}) | < K 
\log (N) / N $$
\end{lem}

\noindent{Proof:}  
Let $Z_i  = 1_{\{X_i \neq Y_i\} }$.  Then $H({X} ) \leq  H({X}, {Z}) =  
H({Y}, {Z} ) \leq  H({Y}) + H({Z} ).$
By symmetry we see that $|  H({X} ) -  H({Y} ) | \leq 
H({Z} ).$ 
But by counting $H({Z}) \leq \log (N^K )$ proving the lemma.   $\Cox$

\noindent{Proof} of Theorem~\ref{entropy}:  
Let $\tilde{B}_n$ have arbitrary boundary conditions and let $B_n$ have 
the 
same vertex set but with unconnected boundary conditions, i.e. the 
equivalence relation consists of singletons.
If $\tilde{\mu}_n$ gives equal probability to each spanning 
forest of $\tilde{B}_n$ in which each component touches the 
boundary then let $\tilde{\nu}_n$
be the measure concentrated on spanning trees of $G$ obtained as follows. 
First partition the vertex set of $G$ with translates 
of $B_n$ and put independent copies of $\tilde{\mu}_n$ on each 
of these translates of $B_n$. Then add (at most $O(n^{d-1})$) edges 
in each translate of $B_n$ to make a tree that spans this translate.  
Then connect each of these trees in the translates 
by a translation-invariant path similar to those constructed in \cite{BK2}.  
This random procedure produces a random spanning tree on $G$ with 
two ends whose measure is denoted $\tilde{\nu}_n$.  It may be made 
translation-invariant by averaging the distribution over all shifts 
by $\Z ^d$ shifts in $B_n$.  By the lemma above we see that for each $n$
$$\frac{1}{|B_n|} H(\tilde{\mu}_n , \tilde{B}_n) \leq H(\tilde{\nu}_n 
) + O (\frac{n^{d-1} 
\log (n) }{n^d}) \leq  H_{top} + O( \frac{n^{d-1} \log (n) }{n^d}).$$
Likewise it is seen that 
$$ \lim \frac{1}{|B_n|} H ( {\mu}_n , {B}_n ) = H_{top} $$ 
where $\mu_n$ is the uniform measure on spanning forests 
of $B_n$ in which every component touches the boundary.  A similar argument 
together with the subadditivity of entropy gives us that $H(\mu_G) = H_{top}$.

Further, if $\mu$ is any ergodic translation-invariant probability with 
$H(\mu) = H_{top}$ then we can show that $\mu = \mu_G$.  Fix a rectangular 
box $B$.  Consider a much larger $C$ around $B$.   
Condition on the $\mu$-outside of $C$ and record the boundary condition
$x \equiv y$ iff $x$ and $y$ are connected by a path in $C^c$.
Since $\mu$ has maximal entropy the conditional distribution of $\mu$ on $B$ 
is the same as the distribution of uniform essential spanning forests 
with these boundary conditions.  (Were this not true we would be able 
to modify $\mu$ within such boundaries and force the entropy to be strictly 
larger.)  If the outside box is large enough then we see (again using 
the arguments
in \cite{Pe} ) that the conditional 
distribution of $\mu$ on $B$  is very close to the distribution of $\mu_G$ 
on $B$ uniformly in the boundary conditions.  Integrating this 
$\mu$-conditional 
distribution with respect to the $\mu$-outside of the large box gives 
us that $\mu$-distribution on $B$ is very close to the $\mu_G$-distribution 
on $B$.  Taking limits gives $\mu = \mu_G$.  This proves part (a) of the 
Theorem and leaves only the computation in (b).

Now we consider finite subgraphs with periodic boundary conditions.  Consider 
the toral graph $(\Z_n)^d  =  \{1,2,...,n\}^d$ with nearest neighbor 
relation taken $\mbox{modulo } n$.  Our vertex set will be 
$(\Z_n)^d \times S$ with incidence matrix
$$M_n ((x,i),(y,j)) =   R^{y-x}(i,j)$$
where the $y-x$ is taken $\mbox{mod }n$ and $n$ is assumed to be large 
enough 
that  
$|x| \geq  n$ implies $R^x =  0.$

Let $\tilde{N}_{B_n} = $  the number of spanning trees on this graph.
We have shown above that $\displaystyle{ H_{top} =  \lim \frac{1}{n^d k} 
\log(\tilde{N}_{B_n})}.$

To complete the proof of the Theorem we use the Matrix Tree Theorem to
compute $\tilde{N}_{B_n}$.  This theorem ( \cite[page 38]{CDS} ) says 
that 
if we have a D-regular connected graph with L vertices and if the eigenvalues 
of the incidence matrix are $\lambda_1, \lambda_2, \ldots , \lambda_{L-1}, 
\lambda_L = D$ then the number of spanning trees of the graph is
$$\frac{1}{L} \prod_{j=1} ^{L-1} (D - \lambda_i).$$

So it is enough to compute the eigenvalues of the matrix $M_n$.  Given $\alpha 
= (\frac{a_1}{n}, \ldots , \frac{a_d}{n})$ for $a_i \in \Z_n$ 
(which we may view as an element of $T^d$) suppose 
that $ \lambda (\alpha)$ is an eigenvalue of $Q( \alpha)$
with eigenvector $v(\alpha)$.  Then $D\lambda(\alpha)$ 
is an eigenvalue of $M_n$ with eigenvector $v(\alpha) \tensor \xi^\alpha$.

This is checked analogously with Lemma~\ref{eigenfunctions}.
\begin{eqnarray*}
&&\sum_{j,y} M_n ((x,i),(y,j)) v( \alpha)_j \exp(2 \pi i \alpha 
  \cdot y) \\
&=&\sum_{j,y} R^{y-x}(i,j) \exp(2 \pi i \alpha \cdot (y-x))  v( \alpha)_j 
  \exp(2 \pi i \alpha \cdot x) \\
&=&\sum_j D Q(\alpha)_{i,j} v(\alpha)_j\exp(2 \pi i \alpha \cdot x) \\
\end{eqnarray*}

Now $D \cdot Q( \alpha)$ is Hermitian, in particular it eigenvectors span 
$\CC^k$ and the characters $\exp (2 \pi i \alpha \cdot x)$ span 
$\CC^{\Z_{n}^{d}}$ 
we see that we have found a complete contingent of eigenvalues.

\begin{eqnarray*}
\tilde{N}_{B_n}&=&\frac{1}{kn^d} \prod_{\lambda, \alpha} (D - D \lambda 
( \alpha)) \\
&=&\frac{1}{kn^d} \left ( \prod_{\alpha \neq 0}  D^k \prod_{\lambda} 
(1 -  \lambda  ( \alpha)) D^{k-1} \right ) \left (
\prod_{\lambda \neq 1} (1 - \lambda(0)) \right ) \\
&=& \frac{1}{kn^d} ( \prod_{\alpha \neq 0} D^k \chi_{Q( \alpha)}(1)) 
D^{k-1} (\prod_{\lambda \neq 1}(1 - \lambda(0))).
\end{eqnarray*}

Continuing and ignoring some logarithmically insignificant terms we get
\begin{eqnarray*}
H_{top}&=&\lim_n \frac{1}{kn^d} \log (\tilde{N}_{B_n}) \\ [2ex]
&=&\lim _n \frac{1}{k} \sum_{\alpha \neq 0} \log (D^k 
\chi_{Q(\alpha)}(1))\frac{1}{n^d} \\ [2ex]
&=&\frac{1}{k} \int_{\VT} \log (D^k \chi_{Q( \alpha)} (1)) d\alpha.
\end{eqnarray*}

The last equality follows from approximating the integral by Riemann 
sums.  Because of the estimate of the eigenvalues from Lemma~\ref{inverse 
square} we see that
$| \chi_{Q(\alpha)}(1)| \leq K |\alpha|^{-2k}$ so that the Lebesgue 
integral is finite.  On $|\alpha| > \ee$ all the terms in the sum
for all $n$ are bounded, so by bounded convergence, the summations
(qua integrals of step functions) converge to the integral.  On
$0 < |\alpha| < \ee$, all sums and the integral go to zero 
as $\ee \rightarrow 0$, implying the desired convergence.    $\Cox$

\section{Dominoes }

A {\em domino tiling} for a graph $G$, otherwise known as a perfect 
matching or a 1-factor, is a collection of edges of $G$ the disjoint
union of whose vertices is $V(G)$, the vertex set of $G$.
There is a correspondence between spanning
trees of a planar graph and domino tilings of a related graph which we
now describe.
Let $G$ be a nice planar graph with vertex set $V_G$ and edge set $E_G$,
``nice'' meaning here that every vertex has finite degree and every
face including the exterior face is bounded by finitely many edges. 
Since $G$ is planar there is a dual graph $G^{\ast}$  with vertices 
$V_{G^\ast}$ and edges $E_{G^{\ast}}$.  In rough terms, $G^{\ast}$  is 
obtained from $G$ by putting a vertex at each face of $G$ and joining 
two such vertices by an edge if their corresponding faces in $G$ meet 
at an edge.  The set $V_{G^\ast}$ is identified with the faces of $G$ 
and the edge sets $E_G$ and $E_{G^{\ast}}$ are also identified.  We construct 
a new bipartite graph $\tilde{G}$ whose vertex set is the union of 
$V_G$, $V_{G^\ast}$ 
and $E_G$.  There is an edge of $\tilde{G}$ joining $v \in V_G$ and 
$e \in E_G$ if and only if $e$ is incident to $v$.
Likewise $v \in V_{G^\ast}$  and $e \in 
E_G$ are joined by an edge if $v$ is a vertex of the edge in $G^{\ast}$  
identified with $e$.  Both $G$ and $G^{\ast}$  sit inside of $\tilde{G}$  
in the sense that $G = \tilde{G} / E_G - V_{G^{\ast}}$ and 
$G^{\ast} = \tilde{G} / E_G - V_G$ (recall this notation from Section 2
for contraction and deletion of a graph).  Here is an illustration of
this where $G$ is the triangular lattice (vertices are filled
circles), $G^*$ is its hexagonal dual (open circles) and the
extra vertices of $\tilde{G}$ are the crossing points. \\
\setlength{\unitlength}{2pt}
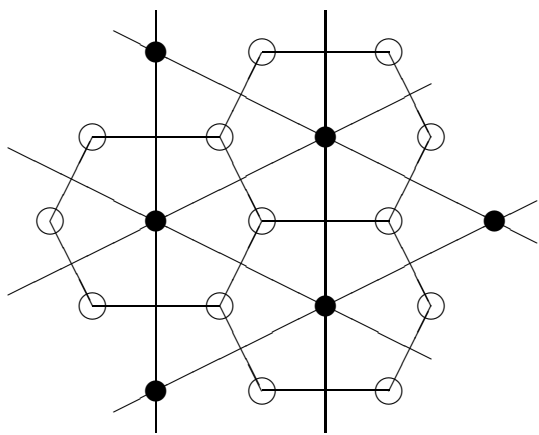
\begin{picture}(160,100)(-40,0)
\put(40,48){\circle {5}}
\put(48,32){\circle {5}}
\put(48,64){\circle {5}}
\put(72,32){\circle {5}}
\put(72,64){\circle {5}}
\put(80,48){\circle {5}}
\put(80,16){\circle {5}}
\put(80,80){\circle {5}}
\put(104,16){\circle {5}}
\put(104,80){\circle {5}}
\put(104,48){\circle {5}}
\put(112,32){\circle {5}}
\put(112,64){\circle {5}}
\put(40,48){\line(1,-2){8}}
\put(40,48){\line(1,2){8}}
\put(72,32){\line(1,-2){8}}
\put(72,32){\line(1,2){8}}
\put(72,64){\line(1,-2){8}}
\put(72,64){\line(1,2){8}}
\put(104,48){\line(1,-2){8}}
\put(104,48){\line(1,2){8}}
\put(104,16){\line(1,2){8}}
\put(104,80){\line(1,-2){8}}
\put(48,32){\line(1,0){24}}
\put(48,64){\line(1,0){24}}
\put(80,16){\line(1,0){24}}
\put(80,48){\line(1,0){24}}
\put(80,80){\line(1,0){24}}
\put(60,16){\circle* {4}}
\put(60,48){\circle* {4}}
\put(60,80){\circle* {4}}
\put(92,32){\circle* {4}}
\put(92,64){\circle* {4}}
\put(124,48){\circle* {4}}
\put(60,8){\line(0,1) {80}}
\put(92,8){\line(0,1) {80}}
\put(52,12){\line(2,1){80}}
\put(52,84){\line(2,-1){80}}
\put(32,34){\line(2,1){80}}
\put(32,62){\line(2,-1){80}}
\end{picture} 

There is a natural correspondance between subgraphs of $G$ and subgraphs 
of $G^{\ast}$ .    If $T$ is any subgraph of $G$ let $T^{\ast}$  be the 
subgraph of $G^{\ast}$  obtained by declaring an edge $e^{\ast}$  
in $T^{\ast}$  
if and only if the corresponding edge $e$ is not in $T$.  Clearly this 
is also a dual operation so $T^{\ast \ast}  = T$.  We record an
easy lemma on dual trees.

\begin{lem} \label{comb lemma}
(a)  Let $G$ be a finite planar graph.  Then $T^{\ast}$  is a spanning 
tree of $G^{\ast}$  if and only if $T$ is a spanning tree of $G$. \\
(b)  Let $G$ be an infinite planar graph all of whose faces are bounded 
regions.  Then $T$ is a one ended spanning tree if and only if $T^{\ast}$  
is a one-ended spanning tree.  Also $T$ is an essential spanning forest 
if and only if $T^{\ast}$  is an essential spanning forest.   $\Cox$
\end{lem}

\noindent{Remark: } This lemma and its soon to be described
connection with domino tilings were noticed independently by Jim Propp 
who suggests the name ``Temperleyan'' for graphs that are isomorphic to  
$\tilde{G}$  for some $G$.  Temperley \cite{Te} first used this trick 
in the case that  $\tilde{G}$  was $\Z^2$.  We also learned of this 
independently from Piet Kastelyn (personal communication).  

Define a {\em directed} essential spanning forest to be a spanning
forest together with a choice of an end for each component.  Think
of edges of a directed spanning forest being oriented toward
this end.  If $T$ and $T^*$ are dual essential spanning forests
of a nice infinite planar graph, say the pair $(T,T^*)$ is directed
if an end has been chosen of each component of $T$ and of $T^*$.
For any nice infinite planar graph $G$, we now describe  
a bijection between domino tilings of $\tilde{G}$ and directed
pairs of essential spanning forests of $G$ and $G^*$.  

If $(T,T^*)$ is such a pair, then let $\Psi (T,T^*)$ be the
domino tiling $A \subseteq E(\tilde{G})$ such that 
\begin{quotation}
$(i)$ the edge from $v \in V(G)$ to $e \in E(G)$ is in $A$ if and only if 
$e \in T$ and is oriented away from $v$, and \\

$(ii)$ the edge from $v^* \in V(G^*)$ to $e \in E(G^*) = E(G)$
is in $A$ if and only if $e \in T^*$ and is oriented away from $v^*$.
\end{quotation}
It is easy to verify that $A$ is a domino tiling: each vertex
$v \in V(G)$ is in precisely one edge of $A$, corresponding to
the unique edge in $T$ out of $v$; similarly each $v^* \in V(G^*)$
is in a unique edge of $A$; and each $e \in E(G)$ is in a unique
edge of $A$ since $e$ is in precisely one of $T,T^*$.  
Conversely, If $A \subset E(\tilde{G})$ is a domino tiling, then 
each edge $f \in A$ connects some $e \in E(G)$ either to
some $v \in V(G)$ or some $v^* \in G^*$.  Let $\Phi (f)$ be the
edge $e$ in either $G$ or $G^*$ accordingly and orient it away
from $v$ or $v^*$.  Then the collection of all $\{ \Phi (f) : f \in A \}$
is the union of a subgraph $G'$ of $G$ and 
the corresponding dual subgraph ${G'}^*$ of $G^*$.
If $G'$ has a loop, then inside the loop is a conponent of ${G'}^*$.
Starting anywhere in this component of ${G"}^*$ and following 
the orientation creates a loop since the component is finite.
This loop encloses a loop of $G'$ inside the original loop.
But this cannot continue forever, whence $G'$ (and ${G'}^*$)
has no loop.  Thus $G'$ and ${G'}^*$ are essential spanning
forests with each component directed toward an end.

Write $\Pi$
for the map that takes a directed pair of ESF's $(T,T^*)$, and forgets
about $T^*$ and about the arrows, producing the undirected ESF $T$.
Then we have established a correspondence
$$\mbox{DOMINO TILINGS} \begin{array}{c} \Phi \\[-2ex] \longrightarrow \\[-2ex]
\longleftarrow \\[-2ex] \Psi \end{array} \mbox{DIRECTED ESF's} \begin{array}
{c} \Pi \\[-2ex] \longrightarrow \\[-2ex] ~ \end{array} \mbox{ESF's} .$$

Now the projection $\Pi$ from directed dual pairs of essential 
spanning forests to essential spanning forests that forgets the
orientation and $T^*$ is not in general one to one, but by Lemma~\ref{comb 
lemma}, if $T$ is a one-ended essential spanning tree of a nice planar
graph, then so is $T^*$, so there is no choice to be made in 
orienting the components of $T$ and $T^*$.   
Now fix a $\Z^2$-periodic planar graph $G$, so the vertex set is 
$\Z^2 \times S$ where $S = \{1, \ldots , k \}$.  There is then
a well defined map $\Psi \circ \Pi^{-1}$ from one-ended spanning trees 
of $G$ to domino tilings of  $\tilde{G}$ (which is also $Z^2$ periodic).  
The uniform spanning forest measure, $\nu$ on $G$ is supported on the 
set of one-ended trees~\cite{Pe} so the above correspondence gives 
a transported measure  $\tilde{\nu}$  on domino tilings of  $\tilde{G}$.

\begin{th} \label{unique max entropy}
The measure  $\tilde{\nu}$ defined above is the unique measure of maximal 
entropy amongst all shift invariant probability measures on domino tilings 
and its entropy per vertex is $kH( \nu )/2e$ where $H( \nu )$ is the 
entropy of $\nu$  and $e$ is the number of edges per fundamental 
domain.
\end{th}

\noindent{Remark:}  This theorem works because boundary conditions are, 
as we have seen in the previous section, irrelevant for trees.  Since 
different boundary conditions for domino tilings give different entropies 
\cite{Elk,Kas,TF}, the excursion through trees is the only soft method 
we know to get uniqueness of the maximal entropy measure for domino
tilings of Temperleyan graphs.  Schmidt \cite[page 58]{Sc} cites an 
argument by Kuperberg that is supposed to prove this for $\Z^2$.

\noindent{Proof:}  We have seen that $\tilde{\nu}$ is well defined and
it is evidently $\Z^2$-invariant, so it remains to prove the assertions
about its entropy.
Suppose that  $\tilde{\mu}$ is a translation invariant probability 
measure on domino tilings of $\tilde{G}$.  This may be transported
to a measure $\mu$ on essential spanning forests of $G$ by 
$\mu (B) = \tilde{\mu}
[\Psi [\Pi^{-1} [B]]]$.  We show that the entropy per fundamental
domain is preserved.  First, note that 
$\mu$  is translation invariant so with probability one the 
components of the essential spanning forest have one or two ends \cite{BK}.  
There is only one way that a one-ended tree may be paved with dominoes 
and there are two ways a two-ended tree may be paved.  Thus the ambiguity 
in determining the domino tiling is one bit for component two-ended tree 
in the forest in $G$ plus one bit for each two-ended component in the 
dual spanning forest of $G^{\ast}$.   Since there are $O (n)$ such components 
in every box of side length $n$ which has on the order of $n^2$ vertices 
we see that the entropy of  $\mu$  and  $\tilde{\mu}$  are the same.  

Now $H(\tilde{\mu})$ per fundamental domain $= H(\mu) \leq H(\nu) = 
H(\tilde{\nu})$ per fundamental domain with equality only when $\mu =
\nu$.  But $\nu$ is concentrated on one-ended spanning trees \cite{Pe}
and hence $\tilde{\nu}$ is the only measure which transports to $\nu$,
which establishes that $\tilde{\nu}$ is the unique measure of maximal
entropy on domino tilings.

Finally, recall $k$ is the number of vertices of $G$ in each fundamental 
domain.  Let $e$ be the number of edges and $f$ the number of faces,
in the sense that a box $\{ 1 , \ldots , n \}^d \times \{ 1 , \ldots , k
\}$ will have approximately $(1+o(1))fn^d$ faces completely contained in it
as $n \rightarrow \infty$.  
Euler's formula applied asymptotically says that $k+f 
= e$.  The entropy of the domino process on  $\tilde{G}$ is the same 
as the entropy of the spanning tree process on $G$ (and as the spanning 
tree process on $G^{\ast}$ ) when measured per fundamental domain.  We 
have given the entropy formula for the spanning tree process per vertex.  
To convert this to the entropy of the domino process on  $\tilde{G}$  
per vertex we must multiply by $k/(k+e+f) = k/(2e)$.   $\Cox$

\begin{cor}
There is a unique measure of maximal entropy on domino tilings of the 
$2$-lattice $Z^2$. and its entropy is $1/4$ the entropy of the spanning 
tree process on $Z^2$.
\end{cor}

\noindent{Remark: }
This entropy number was first calculated by Kastelyn \cite{Kas} in 1961 
as the exponential growth rate of the number of tilings of large rectangle 
or torus.  Since these are atypical boundary conditions this does not 
necessarily prove that this is the largest entropy possible.

\noindent{Proof:  }
In the theorem take $G = Z^2$ so that $G^{\ast}$ is isomorphic to $Z^2$ 
and  $\tilde{G}$  is also isomorphic to $Z^2$.  The uniform spanning tree 
measure $\nu$ on $G$ induces a measure $\tilde{\nu}$ on domino tilings
of $\tilde{G}$ that has 
a fundamental domain of four vertices.  The measure on  $\tilde{G}$  is 
invariant under the induced $\Z^2$ action; since there a four vertices
of $\tilde{G}$ in a fundamental domain of $G$, this is a subgroup
of index $4$ in the usual group of translations of $\tilde{G}$.
Actually, though it must be invariant under all graph automorphisms
$\sigma$ including $90^\circ$ rotations, since otherwise $\tilde{\nu}
\circ \sigma$ would be a measure on domino tilings distinct from
$\tilde{\nu}$ but with the same entropy, violating the uniqueness
shown in Theorem~\ref{unique max entropy}.   $\Cox$

The fact that $\Pi$ is not in general continuous leads
to a problem in trying to compute f.d.m.'s of $\tilde{\nu}$.
If $\C = \{ T : e_1 , \ldots , e_k \in \Tree \}$ is a cylinder 
event in the space of essential spanning forests, then $\Psi [
\Pi^{-1} [\C ]]$ is a finite union of cylinder events in the space of
domino configurations.  On the other hand, 
if $\tilde{\C} = \{ A : e_1 , \ldots , e_k \in A \}$ is a cylinder event 
in the space of domino configurations, then $\tilde{\C}$ is not necessarily
$\Psi [\Pi^{-1} [\C ]]$ for some elementary cylinder event $\C$
on essential spanning forests.  Thus knowledge of the f.d.m.'s of
$\tilde{\nu}$ would yield the f.d.m.'s of $\nu$ quite directly, but
unfortunately not {\em vice versa}.  To illustrate this, let 
$\tilde{\C}$ be the event of finding a square of two vertical 
dominos with the origin at the lower left corner.  

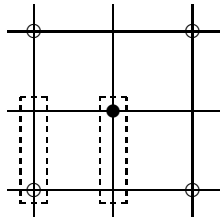
\setlength{\unitlength}{1pt}
\begin{picture}(420,120)
\put (30,40){\line (1,0){80}}
\put (30,70){\line (1,0){80}}
\put (30,100){\line (1,0){80}}
\put (40,30){\line (0,1){80}}
\put (70,30){\line (0,1){80}}
\put (100,30){\line (0,1){80}}
\put (40,40){\circle {5}}
\put (40,100){\circle {5}}
\put (100,40){\circle {5}}
\put (100,100){\circle {5}}
\put (70,70){\circle* {5}}
\put (35,35){\dashbox{2}(10,40)}
\put (65,35){\dashbox{2}(10,40)}
\put(200,80){vertices of $G$ are open circles}
\put(200,60){vertices of $G^*$ are filled circles}
\end{picture}

Then $\Phi [\C]$ is the event that there is an oriented edge from
upward from the origin in $T$ and a dual edge in $T^*$ oriented
downward on the right of the origin.

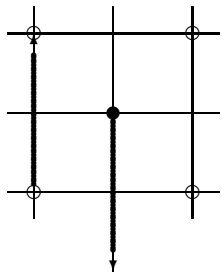
\begin{picture}(120,120)
\put (30,40){\line (1,0){80}}
\put (30,70){\line (1,0){80}}
\put (30,100){\line (1,0){80}}
\put (40,30){\line (0,1){80}}
\put (70,30){\line (0,1){80}}
\put (100,30){\line (0,1){80}}
\put (40,40){\circle {5}}
\put (40,100){\circle {5}}
\put (100,40){\circle {5}}
\put (100,100){\circle {5}}
\put (70,70){\circle* {5}}
\multiput (40,43)(0,1){50}{\circle* {2}}
\multiput (70,67)(0,-1){50}{\circle* {2}}
\put (40,40){\vector(0,1){60}}
\put (70,70){\vector(0,-1){60}}
\end{picture}

The corresponding event on trees is that the edge upward from the origin
be in $T$, that the path from the origin to infinity be through that
edge, that the edge leading right from the origin not be in $T$ and
that the path connecting the origin and the point to the right go
over the top, rather than around the bottom (speaking homotopically
in the plane minus the edge leading right from the origin).  We do not
know how to compute this probability.  There are however some cylinder
domino events corresponding to events whose probabilities we do know how
to compute.  Here is one example.

Consider the following contour which can be broken down into dominos
in the four ways shown. \\
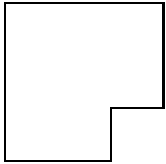
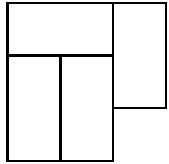
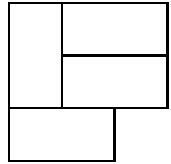
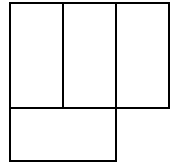
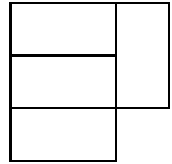
\begin{picture}(500,145)(30,0)
\multiput(20,20)(100,0){5}{\line(1,0){40}}
\multiput(20,20)(100,0){5}{\line(0,1){60}}
\multiput(20,80)(100,0){5}{\line(1,0){60}}
\multiput(80,80)(100,0){5}{\line(0,-1){40}}
\multiput(80,40)(100,0){5}{\line(-1,0){20}}
\multiput(60,40)(100,0){5}{\line(0,-1){20}}
\put(120,60){\line(1,0){40}}
\put(140,60){\line(0,-1){40}}
\put(160,40){\line(0,1){40}}
\put(220,40){\line(1,0){40}}
\put(240,40){\line(0,1){40}}
\put(240,60){\line(1,0){40}}
\put(320,40){\line(1,0){40}}
\put(340,40){\line(0,1){40}}
\put(360,40){\line(0,1){40}}
\put(420,40){\line(1,0){40}}
\put(420,60){\line(1,0){40}}
\put(460,40){\line(0,1){40}}
\put(43,110){$\tilde{A}$}
\put(143,110){$\tilde{A_1}$}
\put(243,110){$\tilde{A_2}$}
\put(343,110){$\tilde{A_3}$}
\put(443,110){$\tilde{A_4}$}
\end{picture}\\
We may calulate the probability of finding this contour with the 
origin at the bottom left.  Since $\tilde{\nu}$ is uniform on the
interior of any box given the boundary, each of the four configurations
inside then has $1/4$ this probability.  To carry this out, map
by $\Phi$ so as to get the following four configurations of directed
edges. \\
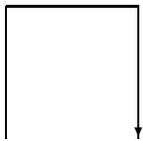
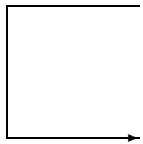
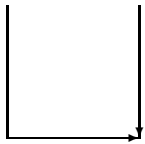
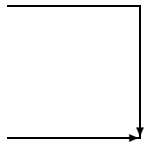
\begin{picture}(400,135)
\put(20,20){\line(0,1){50}} 
\put(20,70){\line(1,0){50}} 
\put(70,70){\vector(0,-1){50}} 
\put(100,20){\vector(1,0){50}} 
\put(100,20) {\line(0,1){50}}
\put(100,70){\line(1,0){50}} 
\put(180,20){\vector(1,0){50}} 
\put(180,20) {\line(0,1){50}}
\put(230,70){\vector(0,-1){50}} 
\put(260,20){\vector(1,0){50}} 
\put(260,70) {\line(1,0){50}}
\put(310,70){\vector(0,-1){50}} 
\put(46,100){$A_1$}
\put(126,100){$A_2$}
\put(206,100){$A_3$}
\put(286,100){$A_4$}
\end{picture} \\
Here, $A_1$ is the event of there being oriented edges in $T$ 
leading up out of the origin, right, and back down, while there 
also being a downwardly directed edge in $T^*$ on the right
of the origin.  Since this dual edge is implied by the the other
three, it is not shown.  A similar thing happens with $A_2, A_3$
and $A_4$.  Now $A_1 \cup A_2 \cup A_3 \cup A_4$ is the event
of $T$ containing three out of the four edges of the square
with the origin at its lower left, and having the path from this
square to infinity exit the square at the lower right.  Then
$\P (\tilde{A}) = \P (A) = \P (A_1 \cup A_2 \cup A_3 \cup A_4)$ which is
by symmetry just $1/4$ times the probability of $T$
containing three of these four edges, with no specified orientation \\
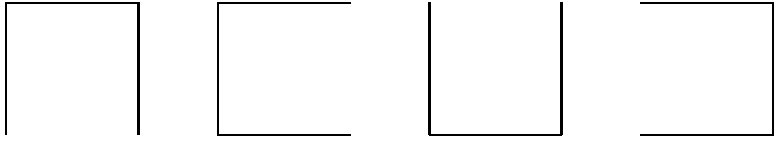
\begin{picture}(400,90)
\put(20,20){\line(0,1){50}} 
\put(20,70){\line(1,0){50}} 
\put(70,70){\line(0,-1){50}} 
\put(100,20){\line(1,0){50}} 
\put(100,20) {\line(0,1){50}}
\put(100,70){\line(1,0){50}} 
\put(180,20){\line(1,0){50}} 
\put(180,20) {\line(0,1){50}}
\put(230,70){\line(0,-1){50}} 
\put(260,20){\line(1,0){50}} 
\put(260,70) {\line(1,0){50}}
\put(310,70){\line(0,-1){50}} 
\end{picture} \\
which is just $1/4$ the probability of the vertex in the center of
the square being a leaf of $T^*$, which is $2\pi^{-2} - 4\pi^{-3}
\approx .0736$ from Section 5. 

The examples which work out this nicely are a small finite class.
There is another class of examples of f.d.m.'s we can calculate.
We give just one illustration, since the taxonomy is still 
being worked out.  Consider the following pair of dominos and
corresponding set of oriented edges of $T$. \\
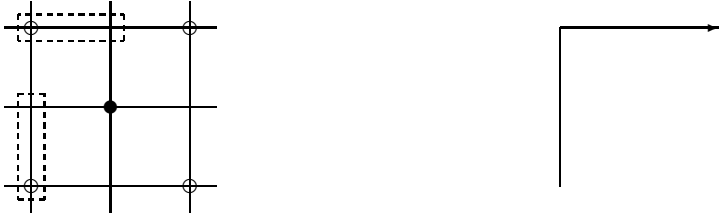
\begin{picture}(420,120)
\put (30,40){\line (1,0){80}}
\put (30,70){\line (1,0){80}}
\put (30,100){\line (1,0){80}}
\put (40,30){\line (0,1){80}}
\put (70,30){\line (0,1){80}}
\put (100,30){\line (0,1){80}}
\put (40,40){\circle {5}}
\put (40,100){\circle {5}}
\put (100,40){\circle {5}}
\put (100,100){\circle {5}}
\put (70,70){\circle* {5}}
\put (35,35){\dashbox{2}(10,40)}
\put (35,95){\dashbox{2}(40,10)}
\put(240,40){\line(0,1){60}}
\put(240,100){\vector(1,0){60}}
\end{picture} \\
The probability of these two oriented edges is the probability
of $T$ containing the two unoriented edges times the conditional 
probability given that of the path to infinity leaving through 
the point at the upper right.  The first probability is computed
by transfer impedances to be
$$\left | \begin{array}{cc} 1/2 & 1/2 - 1/\pi \\ 1/2 - 1/\pi & 1/2
\end{array} \right | \; = \pi^{-1} - \pi^{-2} .$$
The conditional probability may be seen to be the probability
that a random walk coming in from infinity hits the set of three
vertices first at the upper right.  The hitting distribution from
infinity on a set of vertices $x_1 , \ldots , x_k$ is proportional 
to the $k$ entries of $(1 , \ldots , 1) M^{-1}$ where $M$
is the Green's matrix i.e. $M_{ij} = H(x_i , x_j)$ (proof: 
use a last exit decomposition from $\{ x_1 , \ldots , x_k \}$
and then invert the linear relations).  Then
the conditional probability in question is $\pi / (6 \pi - 8)$,
which gives a total probability of $(1 - \pi^{-1})/(6 \pi -8)
\approx .0628$.

\section{Appendix}

\subsection{The classical Green's function}

Let $H (x,y)$ be the classical Green's function for $G$ defined by
\begin{equation} \label{Green transient}
H(x,y) = \sum_{n=0}^\infty \P (SRW_x (n) = y) 
\end{equation}
when $d = 3$, and
\begin{equation} \label{Green recurrent}
H(x,y) = \sum_{n=0}^\infty [\P (SRW_x (n) = y) - \P (SRW_x (n) = x) ] 
\end{equation}
when $d = 1$ or $2$.  
\begin{th} \label{well defined}
The sums in (\ref{Green transient}) and (\ref{Green recurrent}) converge.
Furthermore, $H$ has the following properties:
\begin{quote}
$(i)$ $(H)$ is symmetric; \\
$(ii)$ $H( x , \cdot)$ is harmonic except at $x$, where its excess is 1; \\
$(iii)$ $H$ is bounded if $d \geq 3$; \\
$(iv)$ $H (x,\cdot) - H(y,\cdot)$ is bounded for fixed $x,y$ if $d \leq 2$.
\end{quote}
\end{th}

\noindent{Remark:}  Theorem~\ref{g=H} follows immediately from this
and Corollary~\ref{unique}.  

\noindent{Proof:}  Begin with the observation that $SRW^G$ is transient
if $d \geq 3$ and recurrent if $d \leq 2$; there are many ways, to 
see this, one being to watch $SRW^G$ only at the times when it hits
$\Z^d \times \{ 1 \}$ which is then a symmetric random walk on $\Z^d$
with $\P (x,y)$ having exponential tails.  When $d \geq 3$, the theorem
is now easy to prove.  The sum converges by definition of transience.
Writing $\P (SRW_x (n) = y)$ as the sum of $D^{-n}$ over paths of
length $n$ from $x$ to $y$ shows by path reversal that this is equal
to $\P (SRW_y (n) = x)$ for each $n$, hence $H$ is symmetric.
Boundedness follows from the fact that $H(x,y)  = \P (SRW_x 
\mbox{ hits } y) H(y,y) \leq H(y,y)$, and from the fact that
$H (y,y)$ takes on only $k$ different values.  

Assume now that $d \leq 2$.  Since $SRW^G$ is recurrent there is
a $\sigma$-finite stationary distribution $\mu$, unique up to constant
multiple \cite{IM}.  It is easy to see that this is uniform.
Furthermore, it is well-known \cite{IM} that for any $x,y,z \in G$, 
the ratio of Cesaro averages converges:
\begin{equation} \label{Doeblin}
{1 \over N} \sum_{n=1}^N \P (SRW_z (n) = x) / {1 \over N} \sum_{n=1}^N 
\P (SRW_z (n) = y) \rightarrow \mu (x) / \mu (y) = 1
\end{equation}
as $N \rightarrow \infty$.  Now fix $x,y,z$ and consider the
Markov chain $\{ Z(n) : n \geq 1 \}$ on the space $\{ x, y \}$ gotten 
by looking at $SRW_z$ only when it at $x$ or $y$.  In other words,
$Z(n) = x$ if the $n^{th}$ visit of $SRW_z$ to $\{ x , y \}$ is at $x$
and $Z(n) = y$ otherwise.  The transition matrix for $Z$ is
$\left ( \begin{array} {cc} a & 1-a \\ 1-b & b \end{array} \right )$, where 
$a = \P (SRW_x \mbox{ hits $x$ before } y)$ and $b = 
\P (SRW_y \mbox{ hits $y$ before } x)$.  It follows 
easily from~(\ref{Doeblin}) that the stationary distribution
for $Z$ must be half at $x$ and half at $y$, from which it
follows that $a = b$ in the transition matrix.  

It is easy to calculate 
\begin{equation} \label{Markov Green}
\sum_{n=N}^\infty [\P (Z(n) = x) - \P (Z(n) = y)] =  
   [\P (Z(N) = x) - \P (Z(N) = y)] / (2-2a) . 
\end{equation} 
Now for any positive integers $L < M$, we have 
\begin{eqnarray*}
&& \sum_{n=1}^M [\P (Z(n) = x) - \P (Z(n) = y)] \\[2ex]
& = & \sum_{n=0}^L [\P (SRW_z (n) = x ) - \P (SRW_z (n) = y)] \\
&& + \E \left [ \E (\sum_{n=L+1}^\tau I(SRW_z (n) = x) - I(SRW_z (n) = y)
    \| SRW_z (L+1)) \right ] ,
\end{eqnarray*}
where $\tau$ is the time of the $M^{th}$ visit to $\{ x,y \}$,
and letting $M \rightarrow \infty$ while using~(\ref{Markov Green}) gives
\begin{eqnarray*}
&  & \sum_{n=0}^L [\P (SRW_z (n) = x ) - \P (SRW_z (n) = y)] \\
& = & [2 \P (SRW_z \mbox{ hits $x$ before } y) - 1) / (2 - 2a) \\
&& - [\P (SRW_z (\tau_L) = x) - \P (SRW_z (\tau_L) = y)] / (2-2a) , 
\end{eqnarray*}
where $\tau_L$ is the first time after $L$ that $SRW_z$ hits 
$\{ x , y \}$.  The last term is converging to zero as $L \rightarrow
\infty$, hence letting $z = x$, the sum in~(\ref{Green recurrent}) 
converges.  Moreover, when $z=x$ the sum converges to $1/(2-2a)$,
and having shown that $a=b$ in the transition matrix, we see that
this is symmetric in $x$ and $y$, proving $(i)$.  Along with the
relation $\P (SRW_x (n) = y) = \P (SRW_y (n) = x)$, this also establishes
that $\sum_{n=0}^\infty [\P (SRW_x (n) = x) - \P (SRW_y (n) = y)]
= 0$.  From the fact that $\P (SRW_x (n) = w) \rightarrow 0$
for any $w$ and from the relation 
$$ \sum_{n=0}^N [\P (SRW_x (n) = y) - \P (SRW_x (n) = x) ] =
\dd_x (y) + D^{-1} \sum_{z \sim y} 
\sum_{n=0}^{N-1} [\P (SRW_x (n) = z) - \P (SRW_x (n) = x) ]$$
it now follows that $H(x,\cdot)$ is harmonic except at $x$ and
has excess $1$ at $x$.  

Finally, to check that $H(x,\cdot) - H(y , \cdot)$ is bounded
for fixed $x$ and $y$, use $\sum_{n=0}^\infty [\P (SRW_x (n) = x) - 
\P (SRW_y (n) = y)] = 0$ to conclude that $H(x,z) - H(y,z) = 
\sum_{n=0}^\infty [\P (SRW_z (n) = x) - \P (SRW_z (n) = y)]$.  This
is just $(2 \P (SRW_z \mbox{ hits $x$ before } y) - 1) / (2 - 2a)$, and
the numerator is bounded between $-1$ and $1$, which proves that $H(x, \cdot)
- H(y , \cdot)$ is bounded.    $\Cox$

%
\subsection{Harnack lemmas}

Lemma~\ref{Harnack} is developed in \cite{La} through a series
of theorems beginning with a local central limit theorem for
SRW on $\Z^d$.  We first remark that 
\cite[Theorem 1.2.1]{La} actually holds for the following more general 
random walk. Let $\{ X_n : n \geq 0 \}$ be an irreducible aperiodic
random walk on $\Z^d$ with symmetric transition probabilities
(i.e. $\P(x,x+a) = \P (x,x-a)$) that decay exponentially
(i.e. $\P (x,x+a) = O(e^{-c|a|})$).  Then the characteristic
function for $X_n$ is still given by $\phi (\theta )^n$ where 
$\phi (\theta )$ is real and equal to $1 - <\theta , \theta>
+ O(<\theta , \theta >^2)$ near zero for some positive definite
form $<,>$.  Then the proof of \cite[Theorem 1.2.1]{La} gives
\begin{th}[Local CLT] \label{local CLT}
Under the above assumptions on $X_n$,
there exists a $C > 0$ and a positive definite form $<,>$ for
which $\P_n (x,y) = C n^{-d/2} e^{-<x-y,x-y> / 2n} (1+O(\min (
n^{-1} , <x-y , x-y>)))$.   $\Cox$
\end{th}
Now let $\{ Y_n : n \geq 0 \}$ be a $SRW^G$ started at the point
$(0,1)$.  Write $Y_n = (X_n , Z_n)$, where $X_n \in \Z^d$ and
$Z_n$ is the projected RW on $S$.  It can be shown that $X_n$ and $Z_n$
are exponentially asymptotically independent in the sense that the
joint distribution of $X_n$ and $Z_n$ is within $e^{-cn}$ in total
variation of the product distribution with the correct marginals.
Applying Theorem~\ref{local CLT} to a time change of $X_n$, it can be 
shown that $X_n$ obeys the same local central limit theorem, 
the correction for the time change being smaller than the error bounds 
in the CLT.  This gives
\begin{th}[Local CLT for G] \label{local CLTG}
Let $Y_n$ be a $SRW^G$.  Then
there exists a $C > 0$ and a positive definite form $<,>$ for
which $$\P_n ((x,i),(y,j)) = C n^{-d/2} e^{-<x-y,x-y> / 2n} (1+O(\min (
n^{-1} , <x-y , x-y>))).$$   $\Cox$
\end{th}
This is sufficient to establish part $(i)$ of Lemma~\ref{Harnack}
along the following lines, as pointed out to us by Maury Bramson
(personal communication).  For $x \in B_n$ and $z \in \partial
B_m$, $\nu_x^{B_m} (z)$ is the sum of probabilities of paths 
starting from $x$ and hitting $\partial B_m$ for the first time
at $z$.  Reversing the paths, shows that this is the expected
occupation of $x$ by a SRW starting from $z$ and killed when it
hits $\partial B_m$ again.  

First suppose $d \geq 3$ and fix $\ee > 0$.  Then the local CLT for $G$
allows us to pick $L > n$ large enough so that 
for $w \in \partial B_L$, the occupation measures at $x$ and $y$ 
for $SRW_w$ will be within a factor of $1+\ee$ of each other
for any $x,y \in B_n$.  If $m$ is then chosen large enough, 
the occupation measure for $SRW_w$ at any point in $B_n$
will be at most $1+\ee$ times the occupation measure for
$SRW_w$ killed upon hitting $\partial B_m$.  Now use the Markov
property to write the occupation measure at $x$ for $SRW_z$ killed
upon hitting $\partial B_m$ as a linear combination over $w$ of 
the occupation measure at $x$ of $SRW_w$ killed upon hitting $\partial
B_m$.  This shows the measures at $x$ and $y$ to be within a factor
of $(1+\ee)^2$, and since $\ee$ was arbitrary this establishes
the Harnack principle $(i)$.  

On the other hand, if $d \leq 2$ then $SRW_G$ is recurrent, then
for any $\ee > 0$ and $n$ there is an $m$ large enough so that 
$\P (SRW_x \mbox{ hits $y$ before } \partial B_m) \geq 1 - \ee$ for
all $x,y \in B_n$.  Then $\nu_x^{B_m} \geq (1 - \ee) \nu_y^{B_m}$
for all $x,y \in B_n$, establishing $(i)$.

The remaining parts of the theorem are derived as follows.
To get $(ii)$ from $(i)$, pick $x \in \partial B_n$ and write 
$$\nu_x^{B_m} = \P (SRW_x \mbox{ does not return to } B_n) 
   \rho_x^{B_m B_n} + \P (SRW_x \mbox{ returns to } B_n) \nu' $$
where $\nu'$ is a mixture over $y \in \partial B_n$ of
$\nu_y^{B_m}$, the mixing measure being given by the return
hitting distribution of $SRW_x$ on $\partial B_n$.  Since
$B_n$ is held fixed, $\nu' (\{z\}) / \nu_x^{B_m} (\{z\})$
is converging to $1$ uniformly in $z$ as $m \rightarrow \infty$.
Since $\P (SRW_x \mbox{ returns to } B_n)$ is bounded away from
one (for fixed $B_n$), solving for $\rho_x^{B_m}$ gives
that $\sup_z \rho_x^{B_m} / \nu_x^{B_m} \rightarrow 1$
as $m \rightarrow \infty$.  To get $(iv)$ from $(ii)$, restate
$(ii)$ as saying that the sum, call it $\pi (x,z,n,m)$, of 
$D^{-|\gamma|}$ over paths $\gamma$
from $x \in B_n$ to $z \in B_m$ that avoid $\partial B_n$
and $\partial B_m$ except at the endpoints 
is equal to $(1+o(1)) f(x) g(z,m)$ 
and functions $f,g$ as $m \rightarrow \infty$.  The restatement
of $(iv)$ is easily seen to be identical by time-reversal.
To get $(iii)$ from $(iv)$, just note that $\nu_x^{B_n}$
is a mixture over $y \in \partial (B_m^c)$ of $\rho_y^{B_n B_m}$.

Finally, to get $(v)$ choose $L > n$ so that $B_L$ contains all
the contracted edges.  For $(i)$ and $(ii)$, write $\nu_x^{B_m}$ and
$\rho_x^{B_m B_n}$ as a mixture over $y \in \partial B_L$ of
$\rho_y^{B_m B_L}$ and observe that $\rho_y^{B_m B_L}$ for the 
contracted graph is equal to $\rho_y^{B_m B_L}$ for the uncontracted 
graph, so that all the measures being mixed are identical
up to a factor of $(1+o(1))$.  For $(iii)$ and $(iv)$, write
$\nu_x^{B_n}$ and $\rho_x^{B_n B_m}$ as a mixture over $y \in \partial 
B_L$ of $\nu_y^{B_n}$ and observe that this time it is the 
mixing measures, which are just $\rho_x^{B_L B_m}$ that are all
within a factor of $(1+o(1))$ as $x$ varies with $m \rightarrow 
\infty$. 

\subsection{Convergence of probability measures on trees via moments}

Propositions~\ref{tight means tight} -~\ref{conv criterion} are
adaptations of classical tightness criteria to the setting of
tree-valued random variables.  The development is brief and essentially 
copied from \cite{Du}.  Theorem~\ref{tree moments} is less trivial 
and the proof is given in full detail.

Fix a positive integer $r$ until further notice and restrict attention
to trees of height at most $r$.  Recall the definitions of
$|t|$, $t \wedge r$ and $N(u;t)$ from Section 5.  We say a family of 
probability distributions $\{ \P_n \}$ on such trees is {\em tree-tight} 
if for all $\ee > 0$ there is a $K$ for which $\lim\sup_n \P_n (|\Tree| > K) 
< \ee$.

\begin{pr} \label{tight means tight}
If $\{ \P_n \}$ is tree-tight then every sequence of measures
from $\{ \P_n \}$ has a subsequence that converges in distribution
to a probability measure.
\end{pr} 

\noindent{Proof:}  Since the $\P_n$ laws of $|\Tree|$ are tight
in the usual sense, every sequence has a subsequence $\P_{n_j}$ for which
the laws of $|\Tree|$ converge in distribution to some probability
measure.  For each $k$, the conditional distribution $(\P_{n_j} \|
|\Tree| = k)$ is finitely supported, hence has a subsequence converging
to a probability measure, and diagonalizing over $k$ gives the desired
subsequence.   $\Cox$

\begin{pr} \label{dom moments}
Let $g,h$ be functions from trees of height at most $r$ to the reals.
Suppose that $g>0$ and that $h(t)/g(t) \rightarrow 0$ as $|t| 
\rightarrow \infty$.  Let $\P_n$ be probability measures on these
trees with $\P_n \dconv \P_\infty$ and $\lim\sup_n \E_n g < 
\infty$, where $\E_n$ is expectation with respect to $\P_n$.  Then $\E_n h 
\rightarrow \E_\infty h$.  
\end{pr}

\noindent{Proof:}  $|\E_n h - E_\infty h| \leq  |\E_n h I(|\Tree| > K) - 
\E_\infty h I(|\Tree| > K)| + |\E_n h I(|\Tree| \leq K)| + 
|\E_\infty h I(|\Tree| \leq K)|$.  The last two of these can be
made small uniformly in $n$ by using $|h|/g \rightarrow 0$ and
choosing $K$ large enough.  The first goes to zero for any fixed $K$.
$\Cox$

\begin{pr} \label{coming back}
Suppose $\P_n \dconv \P_\infty$  and for each tree $t$,
$\E_n N(\Tree ; t)$ converges to a finite limit $m(t)$.  Then
$\E_\infty N(\Tree ; t) = m(t)$ for all $t$.  
\end{pr}

\noindent{Proof:}  For any tree $t$ of height at most $r$, let $\C (t)$ 
be the set of all $t'$ of height at most $r$ that extend $t$ by
adding to some vertex in $t$ a single finite chain of descendants.   
Let
$$g_t (u) = \sum_{t' \in \C (t)} N(u ; t') .$$
Notice that $g_t (u) \geq (|u| - |t|) N(u;t)$ since each tree-map 
$\phi : t \rightarrow u$ can be extended to a tree-map of some
$t' \in \C (t)$ into $u$ in at least as many ways as there are
vertices in $u \setminus \mbox{ Image } (\phi)$.  Then
$$\limsup_n \E_n g_t (\Tree) \leq \sum_{t' \in \C (t)} \limsup_n \E_n
   N(\Tree ; t) \leq \sum_{t' \in \C (t)} m(t) < \infty .$$ 
Applying the previous proposition with $g = g_t$ and
$h = N(\cdot ; t)$ finishes the proof.   $\Cox$

\begin{pr} \label{conv criterion}
Suppose $\E_n N(\Tree ; t) \rightarrow \E_\infty N(\Tree ; t) < \infty$
for all $t$ and that $\P_\infty$ is uniquely determined by
the values of $\E_\infty N( \Tree ; t)$.  Then $\P_n \dconv \P_\infty$.
\end{pr}

\noindent{Proof:}  First notice that $\{ \P_n \}$ is tree-tight: letting
$t_j$ be a single chain of $j+1$ vertices, $|\Tree| = \sum_{j=0}^r
N(\Tree ; t_j )$, so $\lim\sup_n \E_n |\Tree| = \lim\sup_n \sum_{j=0}^r
\E N(\Tree ; t_j) < \infty$ by hypothesis, and using $\P (|\Tree| \geq
k) \leq k^{-1} \E |\Tree|$ establishes tightness.  Now each sequence
in $\{ \P_n \}$ has a subsequence converging to a probability
measure $G$ and the previous proposition shows that $\E_n N(\Tree ; t)
\rightarrow \E_G N(\Tree ; t)$ for each $t$.  Then $G = \P_\infty$
by the uniqueness assumption.    $\Cox$

For a positive integer-valued random variable $X$, the moments
of $X$ determine its distribution at least under a condition on the
rate of growth of these moments.  The remainder of the development 
of the method of moments for trees is to prove the following analogous
fact for for tree-valued random variables.  Let $U$ be a random variable 
taking values in the space of locally finite rooted trees of height
at most $r$.  
\begin{th} \label{tree moments}
Suppose $\E N(U ; t)$ is bounded by $e^{k|t|}$ for some $k$.  Then 
the law of $U$ is uniquely determined by the values of
$\E N(U;t)$ as $t$ varies over finite rooted trees of height at most $r$. 
\end{th}
The proof is based on the following version of the integer-valued case.
Let $(A)_s = A (A-1) \cdots (A-s+1)$ denote the $s^{th}$ lower
factorial of $A$.

\begin{lem} \label{high moments}
Suppose $Y \in \Z^+$ and $X$ are random variables and for some
fixed $s > 0$ suppose the values $c_j = \E X(Y)_{s+j}$ exist
for $j \geq 0$ and are bounded above by $e^{kj}$ for some $k$.
Then $\E X I(Y=s)$ is uniquely determined by the values of
the $c_j$'s.  
\end{lem}

\noindent{Proof:}  Let $Z = X(Y)_s$, so that $c_j = \E Z (Y-s) \cdots
(Y-s-j+1)$.  Then the values $d_j \deq \E Z (Y-s)^j$ are determined
from the values of $c_i, i \leq j$ by linear combination.  The
coefficients are bounded by some exponential $e^{kj}$, so the $d_j$ 
are all bounded by some exponential $e^{kj}$.  Let 
$$h(t) = \sum_{j \geq 0} (it)^j d_j / j! = \E Z e^{it(Y-s)} .$$
The power series converges for all $t$, uniformly for $t \in [0,2\pi]$
and hence 
$$\E Z I(Y-s = 0) = (1/2\pi) \int_0^{2\pi} h(t) dt , $$
which yields
$$\E X I(Y=s) = (1/2\pi s!) \int_0^{2\pi} h(t) dt . $$
$\Cox$ 

\noindent{Proof} of Theorem~\ref{tree moments}:  The proof is
by induction on $r$.  The induction will in fact show that
for any random variable $f$ and any $r$, if $U$ has height at 
most $r$ then the values of $\E N(U;t) f$ as $t$ varies over 
trees of height at most
$r$ determine the values of $\E I(U = t) f$, provided that
$\E N(U;t)f$ is bounded by $e^{k|t|}$ for some $k$; using
$f \equiv 1$ will then prove the theorem.  The initial step $r=0$ is 
trivial, since then $U \wedge r$ is always a single vertex.  

Assume then for induction that the theorem is true for some $r$.  
For any tree $t$ let $b(t)$ be the number of children of the root
of $t$.  Define a random vector $Z = 
(Z_1^U , \ldots , Z_{b(U)}^U)$ whose length is always $b(U)$ and 
whose distribution conditional upon $b(U)$ is uniform over
all $b(U)!$ permutations of the $b(U)$ subtrees below the children
of the root of $U$.  In other words $Z$ is this multiset of
subtrees presented in uniform random order.

For any $i$ and $s$ with $s \geq i \geq 0$ and any 
$t_1 , \ldots , t_s$ of height at most $r$, let
$$g_i = g_i (t_1 , \ldots , t_s) = I(Z_1^U = t_1) \cdots
I(Z_i^U = t_i) N(Z_{i+1}^U ; t_{i+1}) \cdots N(Z_s^U ; t_s) I(b(U)=s) .$$
We set up a second induction on $i$ to show that for any $t_i , \ldots ,
t_s$ of height at most $r$ and any random variable $f$, the value
of $\E g_i (t_1 , \ldots , t_s) f$ is uniquely determined by
the values of $\E N(U;t)f$.  

For the initial step $i=0$, choose any $t_1 , \ldots , t_s$ and $f$ 
and write $t_*$ for the tree whose root has $s$ children with 
subtrees $t_1 , \ldots , t_s$.  Write $t_{*,j}$ for a copy of
$t_*$ to which has been added $j$ leaves that are children of the root.
Observe that a map
from $t_{*,j}$ into $U$ is given by choosing $s$ ordered distinct
subtrees $u_1 , \ldots , u_s$ from children of the root of $u$,
mapping each $t_i$ into $u_i$, and then choosing an ordered 
$j$-tuple of vertices from the remaining $b(U)-s$ children of
the root of $U$.  Thus
$$N(U;t_{*,j}) = \sum_{u_1 , \ldots , u_s} N(u_1;t_1) \cdots N(u_s;t_s)
  (b(U)-s)_j $$
where the sum is over ordered $s$-tuples of subtrees from distinct 
children of the root of $U$.  For each $(u_1 , \ldots , u_s)$,
$\P (Z_1 , \ldots , Z_s) = (u_1, \ldots u_s) = 1 / (b(U))_s$. 
Consequently, 
\begin{eqnarray*}
&& \E N(Z_1 ; t_1) \cdots N(Z_s , t_s) (b(U))_{s+j} f \\[2ex]
& = & \E {1 \over (b(U))_s} \sum_{u_1 , \ldots , u_s} N(u_1;t_1) \cdots 
  N(u_s;t_s) (b(U)-s)_j { (b(U))_{s+j} \over (b(U)-s)_j} f \\[2ex]
& = & \E {1 \over (b(U))_s} N(U;t_{*,j}) { (b(U))_{s+j} \over 
    (b(U)-s)_j} f \\[2ex]
& = & \E N(U;t_{*,j}) f .
\end{eqnarray*}
As $j$ varies, $|t_{*,j}|$ increases linearly with $j$. 
By the hypothesis of the theorem this implies
that $\E N(U;t_{*,j}) f$ is bounded by $e^{kj}$ for some $k$.  Then
we may apply Lemma~\ref{high moments} with $Y=b(U)$
and $X= \E N(Z_1 ; t_1) \cdots N(Z_s , t_s) f$ to get that 
these expectations uniquely determine $\E g_0 = N(Z_1 ; t_1) \cdots 
N(Z_s ; t_s) I(b(U)=s) f .$  

Now assume for induction that for any $t_1 , \ldots , t_s$, 
$\E g_i f$ is determined by the values of $\E N(U;t)f$ and write 
$$\E g_{i+1} = \E I(Z_{i+1} = t_{i+1}) [I(Z_1 = t_1) \cdots I(Z_i = t_i) 
   N(Z_{i+2} ; t_{i+2}) \cdots N(Z_s ; t_s) I(b(U)=s) f] . $$
We may apply Lemma~\ref{high moments} with $Y = N(Z_{i+1};t_{i+1})$ 
and $X$ being the rest of the RHS, provided $\E X(Y)_{s+j}$ is
bounded by $e^{kj}$ for some $k$.  But $\E X (Y)_j \leq 
\E N(U;T) f$ where $T$ is a tree whose root has children with 
subtrees: $j$ copies of $t_{i+2}$ and one copy of $t_1 , \ldots t_i$ and
$t_{i+1} , \ldots , t_s$.  Since we have assumed that $\E N(U;t)f$
is bounded by some $e^{k|t|}$, $\E N(U,T) f$ must be bounded by some 
$e^{k'j}$, whence the lemma applies and $\E g_{i+1}$ is indeed determined,
completing the induction on $i$. 

Setting $i=s$ and $f \equiv 1$ now showns that for any $s$ and
any $t_1 , \ldots t_s$, $\P (b(U)=s , Z_1^U = t_1 , \ldots , Z_s^U =
t_s)$ is determined by the values of $\E N(U;t)$.  But this probability
is just $C(t_*) / s!$ times $\P (U = t_*)$, where $C(t_*)$ is the number
of permutations of $\{ 1 , \ldots s \}$ for which $t_i = t_{\pi (i)}$
for all $i$.  Thus $\P (U = t_*)$ is determined for an arbitrary
$t_*$ of height $r+1$, so the induction on $r$ is completed
and the theorem is proved.    $\Cox$

\renewcommand{\baselinestretch}{1.0}\large

\begin{flushright}
Department of Mathematics , Kidder Hall \\
Oregon State University \\
Corvallis, OR 97331-4605 \\
~~\\
Department of Mathematics, Van Vleck Hall \\
University of Wisconsin \\
480 Lincoln Drive \\
Madison, WI 53706
\end{flushright}

\end{document}